\documentclass[runningheads,a4paper]{llncs}

\usepackage{graphicx}
\usepackage{float}
\usepackage{amsmath}
\usepackage{amsfonts}
\usepackage{fontenc}
\usepackage{amssymb,amsmath,array}
\usepackage{todonotes}
\usepackage{textcomp}
\usepackage{caption}

\makeatletter
\renewenvironment{table}%
  {\renewcommand\familydefault\sfdefault
   \@float{table}}
  {\end@float}
\makeatother

\usepackage{url}
\urldef{\mailsa}\path|{alfred.hofmann, ursula.barth, ingrid.haas, frank.holzwarth,|
\urldef{\mailsb}\path|anna.kramer, leonie.kunz, christine.reiss, nicole.sator,|
\urldef{\mailsc}\path|erika.siebert-cole, peter.strasser, lncs}@springer.com|

\begin{document}

\mainmatter  

\title{Anomaly detection for bivariate signals}

\titlerunning{Anomaly detection for bivariate signals}

%
%
\author{Marie Cottrell \and Cynthia Faure \and J\'er\^ome Lacaille \and Madalina Olteanu}
%

\institute{SAMM, EA 4543\\
Panth\'eon-Sorbonne University\\
90 rue de Tolbiac, 75013 Paris, France\\
\url{http://samm.univ-paris1.fr}\\
\&\\
Safran Aircraft Engines, \\
Rond Point Ren\'e Ravaud, R\'eau, 77550 Moissy Cramayel, France\\
\url{https://www.safran-aircraft-engines.com}\\
\&\\
Aosis Consulting, \\
20 impasse Camille Langlade, 31100 Toulouse, France\\
\url{http://www.aosis.net/}
}

%
%

\toctitle{Lecture Notes in Computer Science}
\tocauthor{Authors' Instructions}
\maketitle

\begin{abstract}
The anomaly detection problem for univariate or multivariate time series  is a critical question in many practical applications as industrial processes control, biological measures, engine monitoring,  supervision of all kinds of behavior. In this paper we propose a simple and empirical approach to detect anomalies in the behavior of multivariate time series. The approach is based on the empirical estimation of the conditional quantiles of the data, which provides upper and lower bounds for the confidence tubes. The method is tested on artificial data and its effectiveness has been proven in a real framework such as that of the monitoring of aircraft engines. 
\end{abstract}

\section{Introduction}

Detecting anomalies in univariate and multivariate time series is a critical question in many practical applications, such as fault or damage detection, medical informatics, intrusion or fraud detection, and industrial processes control. The present contribution stems from a joint work with the Health Monitoring Department of Safran Aircraft Engines Company. The motivation behind this collaboration was to find a judicious framework for mining the multivariate high-frequency data recorded on board computers during flights, and isolate unusual patterns, abnormal behaviors of the engine, and possibly anomalies. The issue of anomaly detection on flight data is not new, and some previous results of the joint work with Safran may be found in \cite{Bellas2014} and \cite{Rabenoro2014}.

In a broader context, the literature on anomaly detection, which can also be related to novelty detection, is quite abundant and has been developed for decades in various fields: machine learning, statistics, signal processing. The techniques used for addressing the matter may use supervised or unsupervised learning, model-based algorithms, information theory or spectral decomposition. For quite an exhaustive review, the reader may refer to \cite{Chandola2009}. 

More particularly, we focus here on the issue of detecting collective anomalies or discords -- an unusual subsequence of a time-series, in contrast to local anomalies which consist in unique abnormal time-instants -- in a multivariate context. Detecting collective anomalies or unusual patterns in univariate time-series has been extensively studied, and we may cite, for instance, algorithms based on piecewise aggregated approximation \cite{Lin2003}, nearest-neighbor distances \cite{Keogh2005}, Fourier or wavelet transforms \cite{chen2008multi}, \cite{Michael00}, Kalman filters \cite{Leith}, ... In the multivariate case, anomaly detection has to take into account both the multivariate aspect of the data, and the temporal span, the possibly existing correlations and dependencies.  Whereas the initial approaches used time series projection \cite{Galeano} and independent component analysis \cite{Baragona} to convert the multivariate time series into a univariate one, or performed separate anomaly detection for each variable \cite{Lakhina}, global approaches have been developed only recently. Among these recent works, one may cite, for instance \cite{Cheng}, who use a kernel-based method for capturing dependencies among variables in the time series, and  \cite{Malhotra2015} who use neural networks for isolating anomalous regions in a multivariate time-series.


This paper adresses the issue of detecting unusual patterns in a multivariate time series context. Unlike some of the cited literature above, we suppose that the data comprises a set of multivariate time-series, which have already been segmented into patterns of unequal lengths, using some change-point detection technique. Our goal is to find the most unusual patterns among them, and for doing so we take the unsupervised learning approach  \textit{(as defined by the AI)}. No hypothesis whatsoever is made on some underlying model, the only constraint is to suppose that one component of the time series, called \textit{key variable} in the sequel, exists, may be distinguished, and its behavior strongly influences the behavior of the rest. The approach we introduce here may be briefly described as follows: first, the initial patterns of the \textit{key variable} are summarized by a fixed number of numerical features, second, they are clustered into an optimal number of clusters using self-organizing maps and the computed numerical features, third, the multi-variate patterns are realigned and synchronized within each cluster, fourth, unusual patterns are extracted after computing confidence tubes from empirical quantiles in each cluster. This global approach has already been introduced in \cite{Faure2018}, and the contribution in the present manuscript relies in the use of conditional first order quantiles for computing the confidence tubes, instead of quantiles computed at each time instant, independently on the past. As will be illustrated in the Experiments section, this conditional approach greatly improves the ratio of false positives detection on simulated examples. 

For simplicity purposes, the bivariate case only is presented here, but the algorithm can be easily extended to higher dimensional data. 

The rest of the paper is organized as follows : Section 2 describes the main steps of the proposed methodology, Section 3 contains results on simulated examples and a comparison between the previous version of the method and the modified one based on conditional quantiles, while Section 4 illustrates the method on real-life dataset stemming from flight data. 

\section{Methodology}

Let $S_a = (X_a, Y_a)$ be a set of bivariate $\mathbb {R}^2$-valued time series, for $a= 1, \ldots A$. For each $1 \leq a \leq A$, the lengths of $X_a$ and $Y_a$ are equal and denoted by $l_a$. Note that the lengths $l_a$ can be different from one time series to another one.

We assume that one of both variables is a \textit{key variable} (easier to observe, with a limited number of different behaviors, which influences the behavior of the other). This hypothesis is not very restrictive since in many processes, there is a measure which gives a first information on all the others variables (water temperature, blood composition, temperature of the core, etc. according to the application field.)

So we will define two successive levels of analysis: the first one deals with the key variable (say $X_a$) and the second one will further take into account the second variable ($Y_a$).

A first idea could be to bring together all the signals $X_a$, to compute an "average signal" and to affirm that all the signals which are far from the mean signal are anomalies. The difficulty arises from the fact that the lengths $l_a$ are different, so how to compute the "average signal"?

Secondly, even if we found a method to define an "average signal", there is no reason for that this one will correctly represent all the signal behaviors. A previous clustering of the signals $X_a$ seems to be mandatory in order to get a limited number of homogeneous clusters within which a representative curve can be defined. 

\subsection{Clustering} \label{clustering}

To overcome the difficulty of dealing with   signals of different lengths, as in Faure et al. 2017-a \cite{Faure2017_a}, we replace each signal $X_a$ by a fixed-length vector of some well-chosen numerical features (as length, mid point, median, variance, variance of first half, variance of second half, mean of first half, mean of second half, etc...). Let $M$ be the number of these numerical features. Any classical clustering algorithm can now be used for these $M$-vectors and to build clusters denoted by $C_1, C_2, \ldots, C_I$, where $I$ is the number of clusters.

Then we propose to use a SOM algorithm which builds the clusters into two steps : first we get a large number of Kohonen classes displayed into a Kohonen map, then there are grouped into a limited number of superclasses named clusters. The number of these clusters is adjusted using an empirical criterion based on the percentage of explained variance, computed by a Hierarchical Ascending Classification of the code-vectors and its dendrogramm.

In fact it could also be possible to avoid the replacement of   the signals by a fixed number of features by using the relational SOM algorithm as in Olteanu et al. 2015 \cite{Olteanu2015}.  It needs to define the dissimilarities between  different length signals by Eq. \eqref{eq:dissimilarity}, as explained in section 2.2. We experimented it and got very similar results.

Each cluster $C_i$ contains a set of signals $X$, denoted by $X_{a}^i$ for simplicity. They are similar from the point of view of their numerical features,  but can have different lengths. 

The idea is then to define a Reference Curve in each of these homogeneous sets, to "represent" the cluster and improve the visualization by coming back to the initial representation of the signals. 

\subsection{Dissimilarity between unequal length curves - Reference Curve}

Firstly let us  define the dissimilarity between two curves with different lengths as in Faure et al. 2017-b \cite{Faure2017_b}. Let $X_{a_{1}}$ and $X_{a_{2}}$ be two curves with lengths $l_1$ and $l_2$ and $l_1 < l_2$.

If $X_{a_1}=(x_1, x_2, \ldots, x_{l_1})$,  one defines its extended version  $\tilde{X}_{a_1}$ with length $l_2 + l_1 + l_2$ as follows

$$\tilde{X}_{a_1}=(\underbrace{x_1, x_1, \ldots, x_1}_{l_2} \vert \underbrace{x_1, x_2, \ldots, x_{l_1}}_{l_1}\vert 
\underbrace{x_{l_1},x_{l_1}, \ldots,x_{l_1}}_{l_2}). $$

Note that if $X_{a_{1}}$ is extracted from a time series, the extension at left and at right will be done by using the previous and next values in the complete series.

Then a $l_2$-points sliding window is defined over 
$\tilde{X}_{a_1}$ to cut $l_2$-points segments: 
$I_q( \tilde{X}_{a_1})=\tilde{X}_{a_1}[q, q+l_2 -1]$ is the section of $\tilde{X}_{a_1}$ taken between indexes $q$ and 
$q+ l_2 -1$, for $q=1, \ldots, l_1 + l_2 +1$.

The dissimilarity between $X_{a_1}$ and $X_{a_2}$ is finally defined as 
\begin{equation}\label{eq:dissimilarity}
diss(X_{a_1},X_{a_ 2})= \min_{q \in {1, \ldots, l_1+l_2-1}} 
\frac{\|I_q( \tilde{X}_{a_1})-X_{a_2}\|}{2 l_2}
\end{equation}

As explained in \cite{Faure2017_b}, this reference curve is that one which minimizes the sum of all the dissimilarities between it and all the other curves of the cluster. Let us denote by $RC_i$ the reference curve of cluster $C_i$, and by $L_i$ its length.

Once the reference curve $RC_i$ is computed, the question is to harmonize the lengths of all the curves of cluster $C_i$.
This step is achieved by  applying to each one a transformation which is a combination of translation, completion and truncation  to get the best proximity between each curve and the reference one.

\subsection{Realignment, synchronization-transformation}

The method is presented in \cite{Faure2017_b}. Let us consider a curve  $X_a^i$ in $C_i$ with  length $l_a^i$. As before, the curve is extended at left by $L_i$ constant values equal to its first value $X_a^i(1)$ and at right by $L_i$ constant values equal to its last value $X_a^i(l_a^i)$. This new object is denoted by 
$\hat{X}_a^i$. 

A $L_i$-points  window sliding over $\hat{X}_a^i$ cuts  segments of length $L_i$ of the extended curve $\hat{X}_a^i)$, denoted by
\begin{equation}\label{eq:sliding}
I_q(\hat{X}_a^i)=\hat{X}_a^i[q, q+L_i -1]
\end{equation}

which is the section of 
$\hat{X}_a^i$ taken between indexes $q$ and 
$q+ L_i -1$ , for $q=1, \ldots, l_a^i + L_i+1$. 

Then 
\begin{equation}\label{eq:diss}
diss(X_a^i, RC_i)=\min_{q \in {1, \ldots, l_a^i + L_i+1}}
\frac{\|I_q( \hat{X}_a^i)-RC_i\|}{2 L_i} 
\end{equation}
and 
\begin{equation}\label{eq:argdiss}
q_a^i=\arg \min_{q \in {1, \ldots, l_a^i + L_i+1}}
\frac{\|I_q( \hat{X}_a^i)-RC_i\|}{2 L_i}.
\end{equation}

%
%

Then the curve $X_a^i$ is replaced by  $I_{q_a^i}(\hat{X}_a^i)$ which is denoted by $\breve{X}_a^i$. This transformation is applied to all the curves of $C_i$ and realigns them with the reference curve $RC_i$. Finally, all the curves $\breve{X}_a^i$ for a given $i$ have the same length $L_i$.

This \textit{synchronization-transformation} is applied to the second signal $Y_a^i$ : it is extended, translated, cut at the same indexes as $X_a^i$. The transformed signal is denoted by 
$\breve{Y_a^i}$ and has also the same length $L_i$.

For the sake of simplicity, we denote by $C_i$ the set of signals 
$X_a^i$ as well as that of the transformed signals 
$\breve{X}_a^i$. The corresponding second components $Y_a^i$ or 
$\breve{Y}_a^i$ define a set $D_i$. We denote by $E_i$ the set of all the couples of transformed signals $\breve{X}_a^i, \breve{Y}_a^i$.

\subsection{Anomaly detection}

The goal is now to detect eventual non typical curves in 
$E_i$. The anomalies can be related to the first component, to the second one, or to both.

The matter is to build  confidence tubes in each set $C_i$ and $D_i$ for a given level of confidence. As the data distribution is unknown, the simplest way is to consider all the values for each value of time $t$ and to build empirical point-by-point confidence tubes. However it could be preferable to keep the information given by the time series structure and to consider conditional quantiles, following Samanta et al. 1989 \cite{Samanta1989}, Charler et al. 2014 \cite{Charler2014}, Charlier et al. 2015 \cite{Charlier2015}.

\subsubsection{Empirical point-by-point confidence tubes - CT detection method} \label{sec::CT}

For each $t$ such that $1 \leq t \leq L_i$, we consider the sets of values $\breve{X}_a^i(t)$ and $\breve{Y}_a^i(t)$. 

For a given level $\alpha$ (typically $\alpha = 5 \%$), 
one denotes by $q^X_{t,\frac{\alpha}{2}} $ 
(resp. $q^Y_{t,\frac{\alpha}{2}}$)  and   
$q^X_{t,1-\frac{\alpha}{2}}$  
(resp. $q^Y_{t,1-\frac{\alpha}{2}}$) the $\alpha$-quantiles  of these sets.

The tube of confidence  $1-\alpha$ of the curves $\breve{X}_a^i$ in $C_i$ is defined as the sets of all the curves entirely comprised between
\begin{itemize}
\item the upper bound defined by 
$q^X_{t,1-\frac{\alpha}{2}}, 1 \leq t \leq L_i$
\item the lower bound defined by 
$q^X_{t,\frac{\alpha}{2}}, 1 \leq t \leq L_i$
\end{itemize}

In the same way, we build the confidence tubes for the curves 
$\breve{Y}_a^i(t)$ in $D_i$.\\

\textbf{Definition} : A curve in $C_i$ or in $D_i$ is an anomaly if and only if it has at least $P$ \% consecutive points out of the corresponding confidence tube.\\ 

The threshold $P$ \% is fixed by the user, and is typically equal to 10\%.

\subsubsection{Conditional quantiles - CQ detection method} \label{sec::CQ}

The idea is to take into account that the set of data points 
$(t,\breve{X}_a^i(t)$ (resp. $(t,\breve{Y}_a^i(t)$ is structured since they are observations of time series.

Let us deal with the first components $\breve{X}_a^i(t)$ 
in $C_i$. The same will apply to the second components of the bivariate signals.

The goal is to estimate the confidence intervals of 
$\breve{X}_a^i(t)$ conditionally to the past of the time series, for  instance $\breve{X}_a^i(t-1)$ .

The number of delays to consider can be determined after analyzing the times series with an auto-regressive model. The proposed method can be easily extended to a deeper past, $t-2, t-3, \ldots$, etc. 

For a given level $\alpha$ (typically $\alpha = 5 \%$), and if we only use  the immediate past, 
we define the conditional quantiles 
$q^X_{t,\frac{\alpha}{2}}(x) $ and 
$q^X_{t,1-\frac{\alpha}{2}}(x)$ 
such that

\begin{itemize}
\item $\mathbb{P} \left(\breve{X}_a^i(t) \leq q^X_{t,\frac{\alpha}{2}}(x) / 
\breve{X}_a^i(t-1)=x\right) = \frac{\alpha}{2}$

\item $\mathbb{P} \left(\breve{X}_a^i(t) \geq q^X_{t,1-\frac{\alpha}{2}}(x) / 
\breve{X}_a^i(t-1)=x\right) = \frac{\alpha}{2}$
\end{itemize}

As the conditional distribution of $\breve{X}_a^i(t)$ given 
$\breve{X}_a^i(t-1)$ is unknown, we have to estimate these quantiles. So we need to discretize the set of values 
${\breve{X}_a^i(t-1), 1 < t\leq L_i}$ to get a sufficient number of values of $\breve{X}_a^i(t)$ for each value of 
$\breve{X}_a^i(t-1)$. 

Once the quantiles are computed for each discretized value 
$^d\breve{X}_a^i(t-1)$ of $\breve{X}_a^i(t-1)$, 
\textbf{a curve in $C_i$ is detected as an anomaly} if and only if the number of couples $(\breve{X}_a^i(t-1), \breve{X}_a^i(t))$ such that 
$$\breve{X}_a^i(t) \notin 
\left[q^X_{t,\frac{\alpha}{2}}(^d\breve{X}_a^i(t-1), q^X_{t,1-\frac{\alpha}{2}}(^d\breve{X}_a^i(t-1)\right]$$
is greater than a certain threshold, fixed by the user (typically 10 \%).

\section{Example of simulated signals}

\begin{figure}[H] 
\begin{minipage}[b]{0.3\textwidth}
\centering
\includegraphics[width=0.95\linewidth]{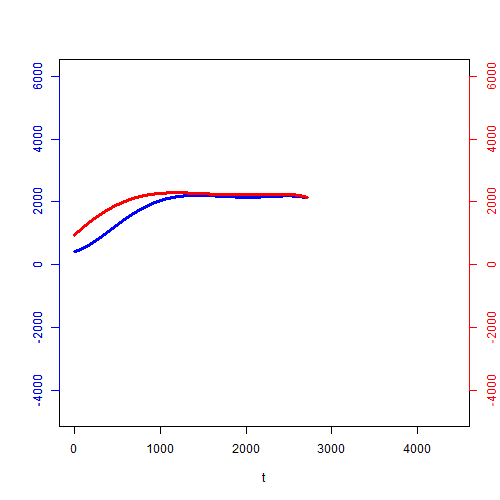} 
\caption{Exemple of bivariate signal, $X$ blue, $Y$ red} 
\label{Fig::sim_curves}
\end{minipage} 
\begin{minipage}[b]{0.3\textwidth}
\centering
\includegraphics[width=0.95\linewidth]{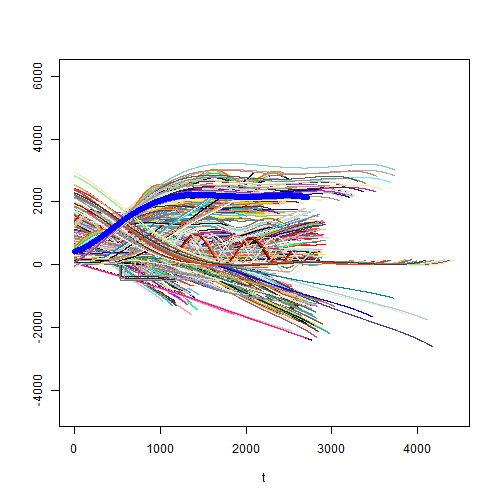} \\
\caption{All signals $X$} 
\label{Fig::sim_curves1}
\end{minipage}
\begin{minipage}[b]{0.3\textwidth}
\centering
\includegraphics[width=0.95\linewidth]{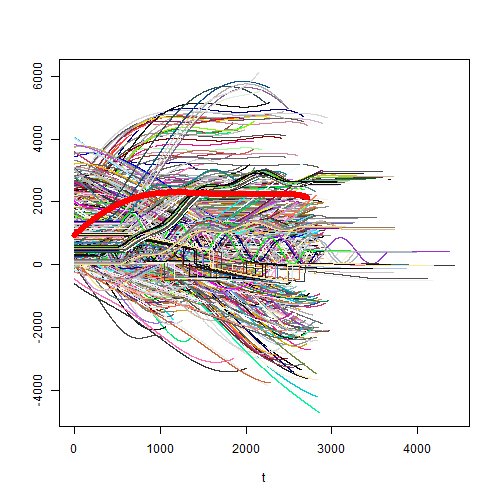}\\ 
\caption{All signals $Y$} 
\label{Fig::sim_curves2}
\end{minipage} 
\end{figure}

We build 2000 artificial bivariate signals shown in Figures \ref{Fig::sim_curves}, \ref{Fig::sim_curves1}, \ref{Fig::sim_curves2}.

\begin{figure}[H]
\centering
\includegraphics[scale=0.45]{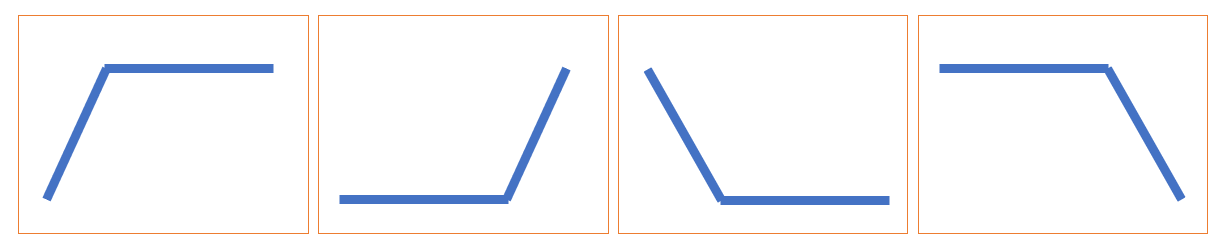}
\caption{Four shapes used for the simulated $X$ signals}
\label{Fig::patterns}
\end{figure}

The $X$ variable has   one of four different shapes which are  shown in Figure \ref{Fig::patterns}.

These shapes are chosen since in the real date framework of flight data, most of the time series can be described as piecewise linear curves.

For each pattern, we randomly choose the length between 500 and 3000 points, the slope change point (between the third and the two third of the length), the slope values which vary from $0.5$ to $4$ (for ascendant shapes) and from $-4$ to $-0.5$ (for descendant shapes). We add a Gaussian centered noise with variance between $10$ and $100$, and we smooth the signal by using a 5 degrees polynomial. See Figure \ref{Fig::sim_curves} some of the simulated curves.

To simulate the second variable $Y$ associated to some variable $X$, we start from  $X$, we randomly choose a breakpoint $z$ and two slopes for the segments before and after this breakpoint, more or less in the same range as for the $x$-curve. If necessary we extend or cut the time series $Y$ to give it the same length as that of X. The values are smoothed by using a 5 degrees polynomial.

\begin{figure}[H] 
\begin{minipage}[b]{0.24\textwidth}
\centering
\includegraphics[width=0.95\linewidth]{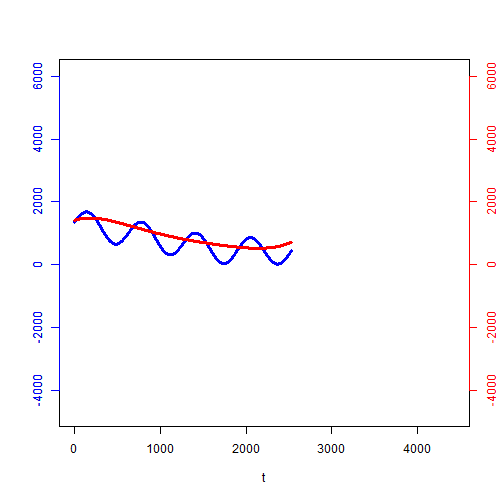} \\
(a)
\end{minipage} 
\begin{minipage}[b]{0.24\textwidth}
\centering
\includegraphics[width=0.95\linewidth]{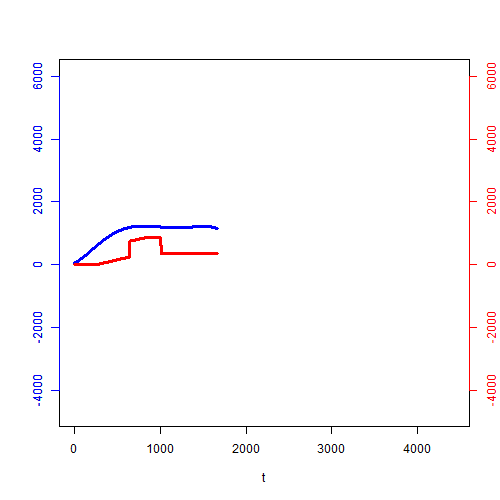} \\
(b)
\end{minipage}
\begin{minipage}[b]{0.24\textwidth}
\centering
\includegraphics[width=0.95\linewidth]{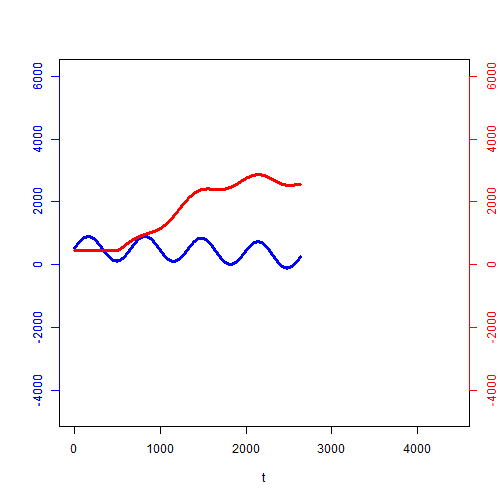} \\
(c)
\end{minipage} 
\begin{minipage}[b]{0.24\textwidth}
\centering
\includegraphics[width=0.95\linewidth]{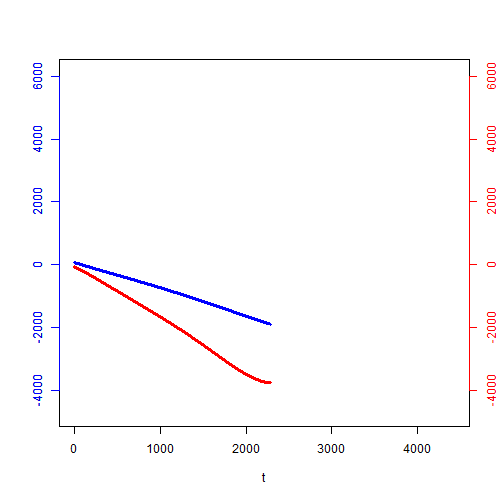} \\
(d)
\end{minipage} 
\caption{Examples of atypical curves ($X$  blue, $Y$ red)}
\label{Fig::weird}
\end{figure}

Among the artificial signals, we simulated 50 atypical shapes (anomalies). We used 4 types of "anomaly shapes": sinusoidal (see Figures \ref{Fig::weird} (a) and (c)), "hat" (see Figure \ref{Fig::weird} (b)) and linear (see Figure \ref{Fig::weird} (d)). A couple of variable $(X,Y)$ can be atypical for $X$ only, for $Y$ only or for both.\\

%

\begin{figure}[H] 
\begin{minipage}[b]{0.19\textwidth}
\centering
\includegraphics[width=0.95\linewidth]{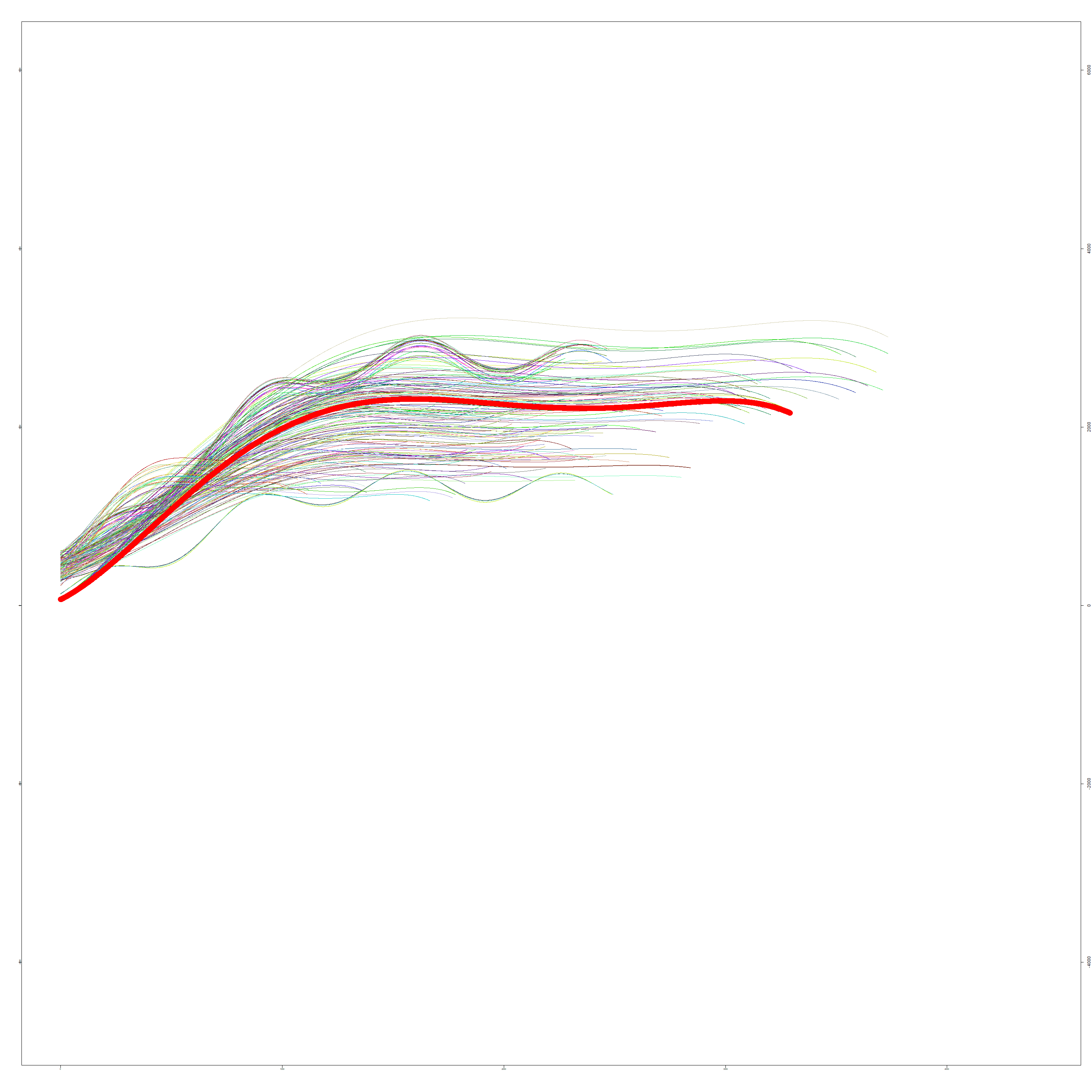} \\
$C_1$
\end{minipage} 
\begin{minipage}[b]{0.19\textwidth}
\centering
\includegraphics[width=0.95\linewidth]{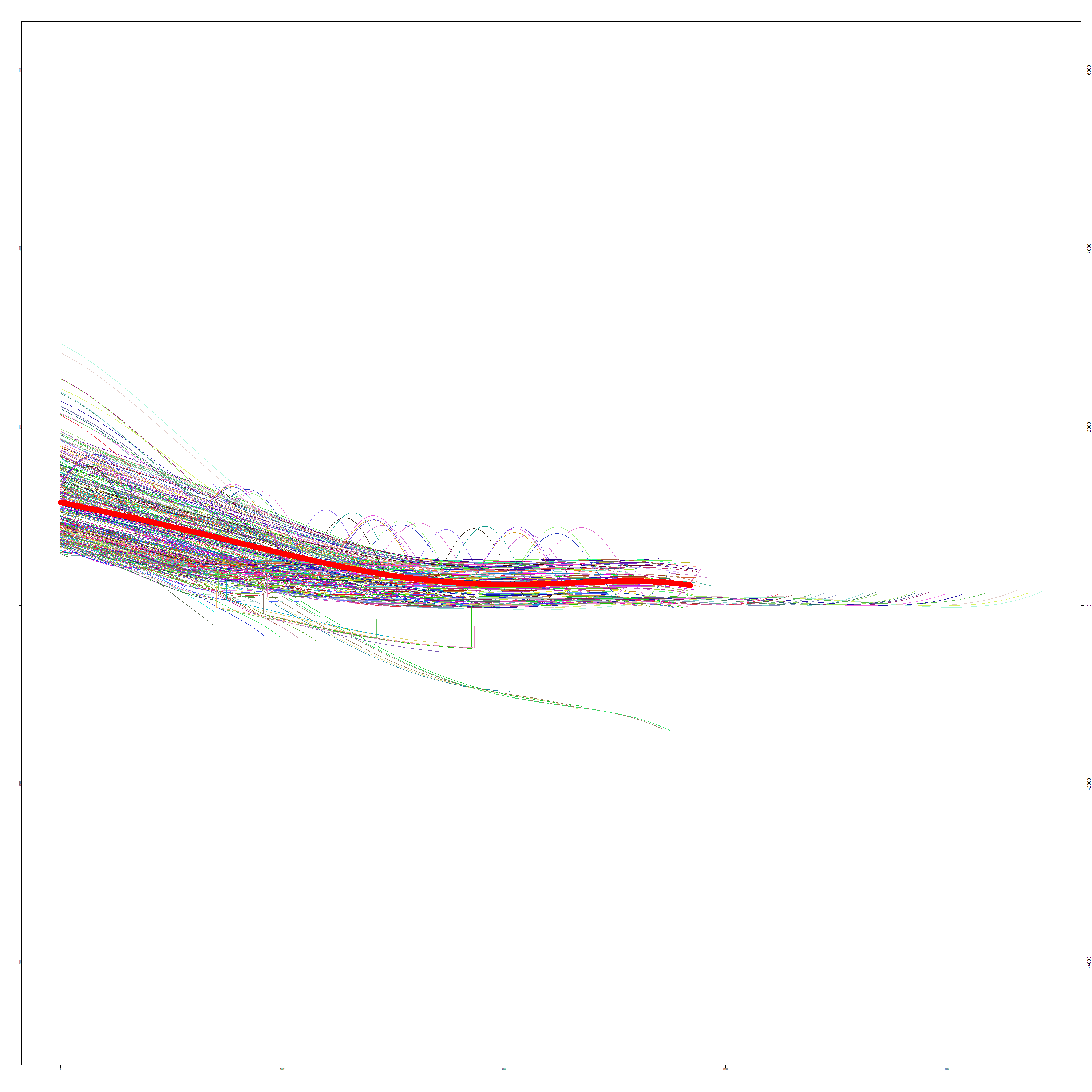} \\
$C_2$
\end{minipage}
\begin{minipage}[b]{0.19\textwidth}
\centering
\includegraphics[width=0.95\linewidth]{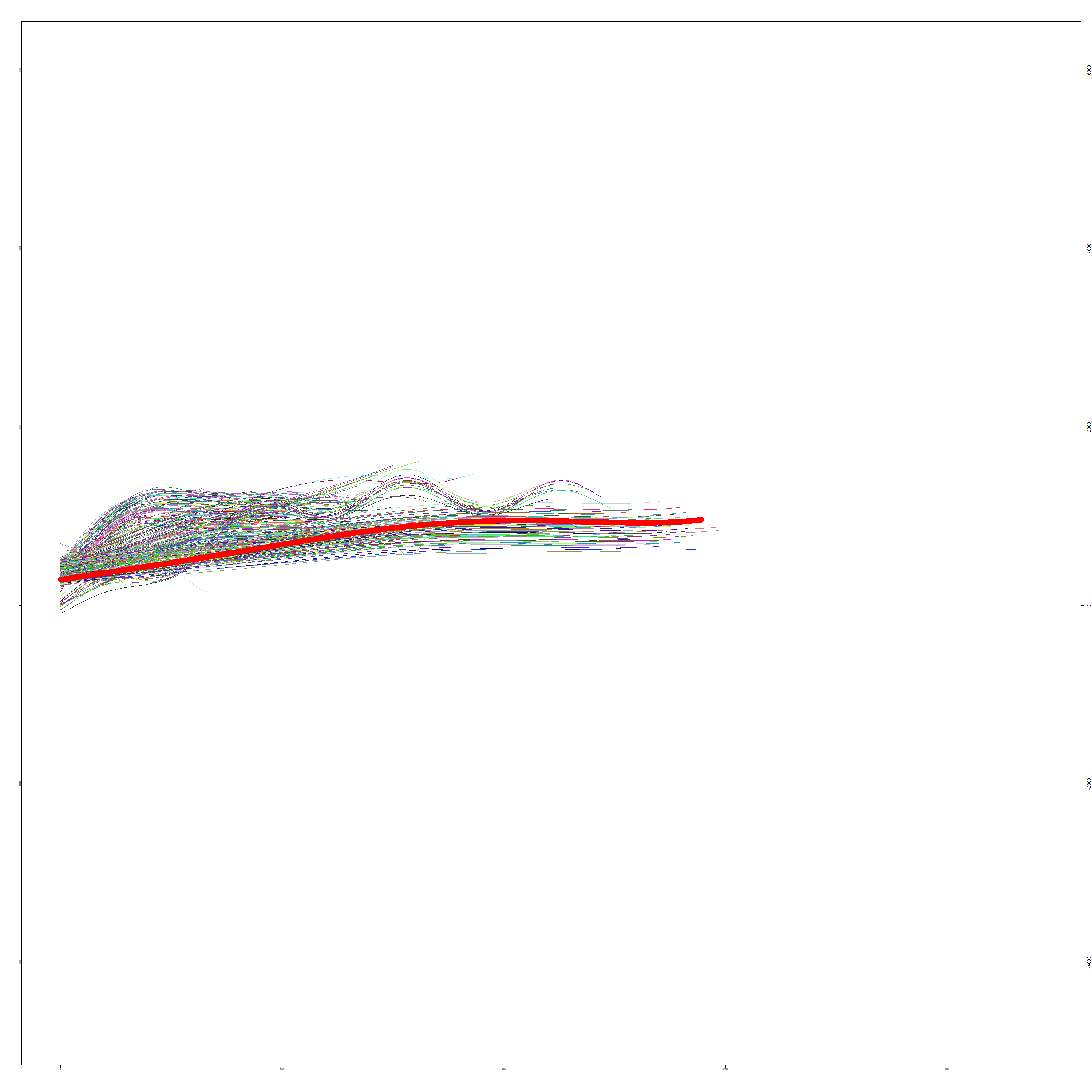} \\
$C_3$
\end{minipage} 
\begin{minipage}[b]{0.19\textwidth}
\centering
\includegraphics[width=0.95\linewidth]{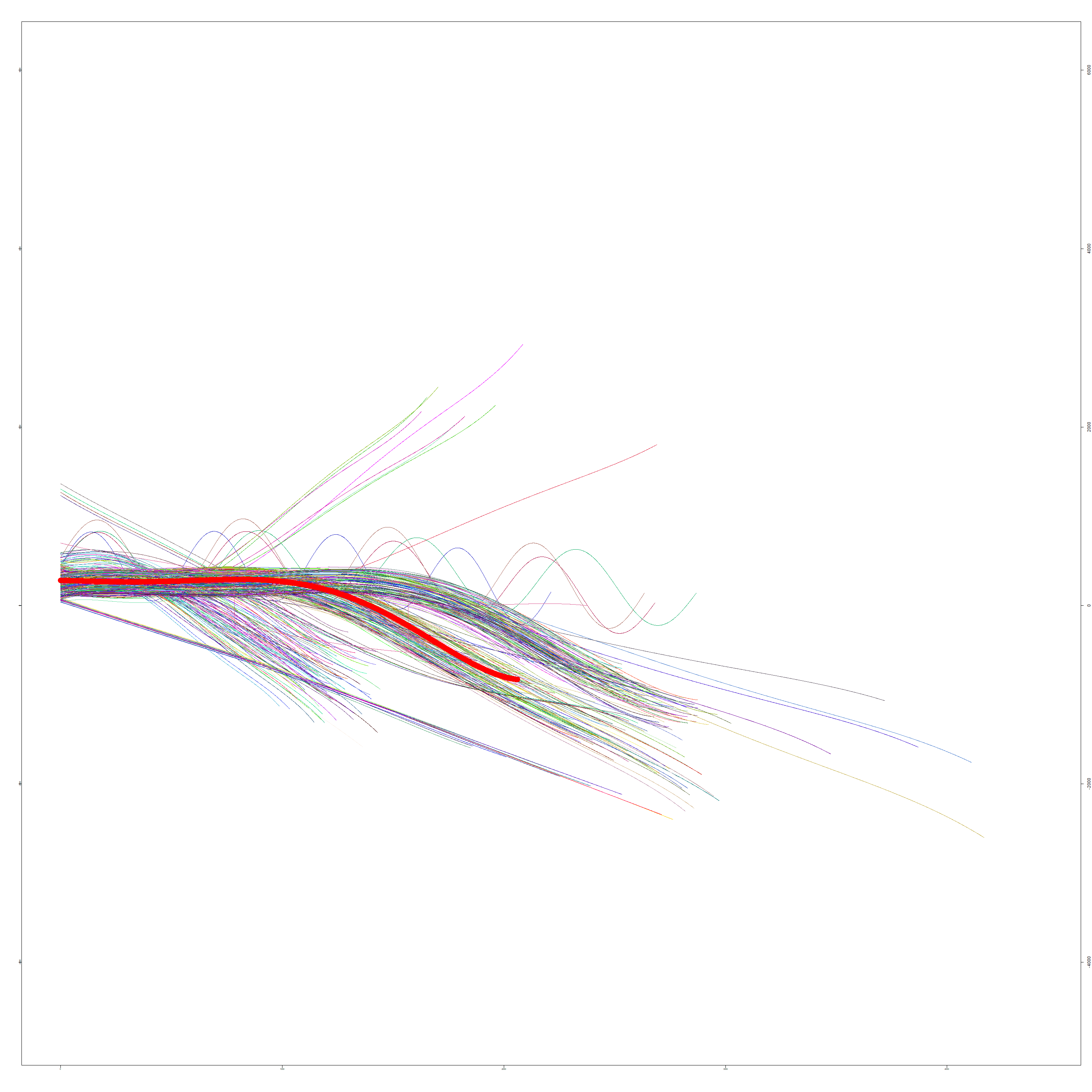} \\
$C_4$
\end{minipage} 
\begin{minipage}[b]{0.19\textwidth}
\centering
\includegraphics[width=0.95\linewidth]{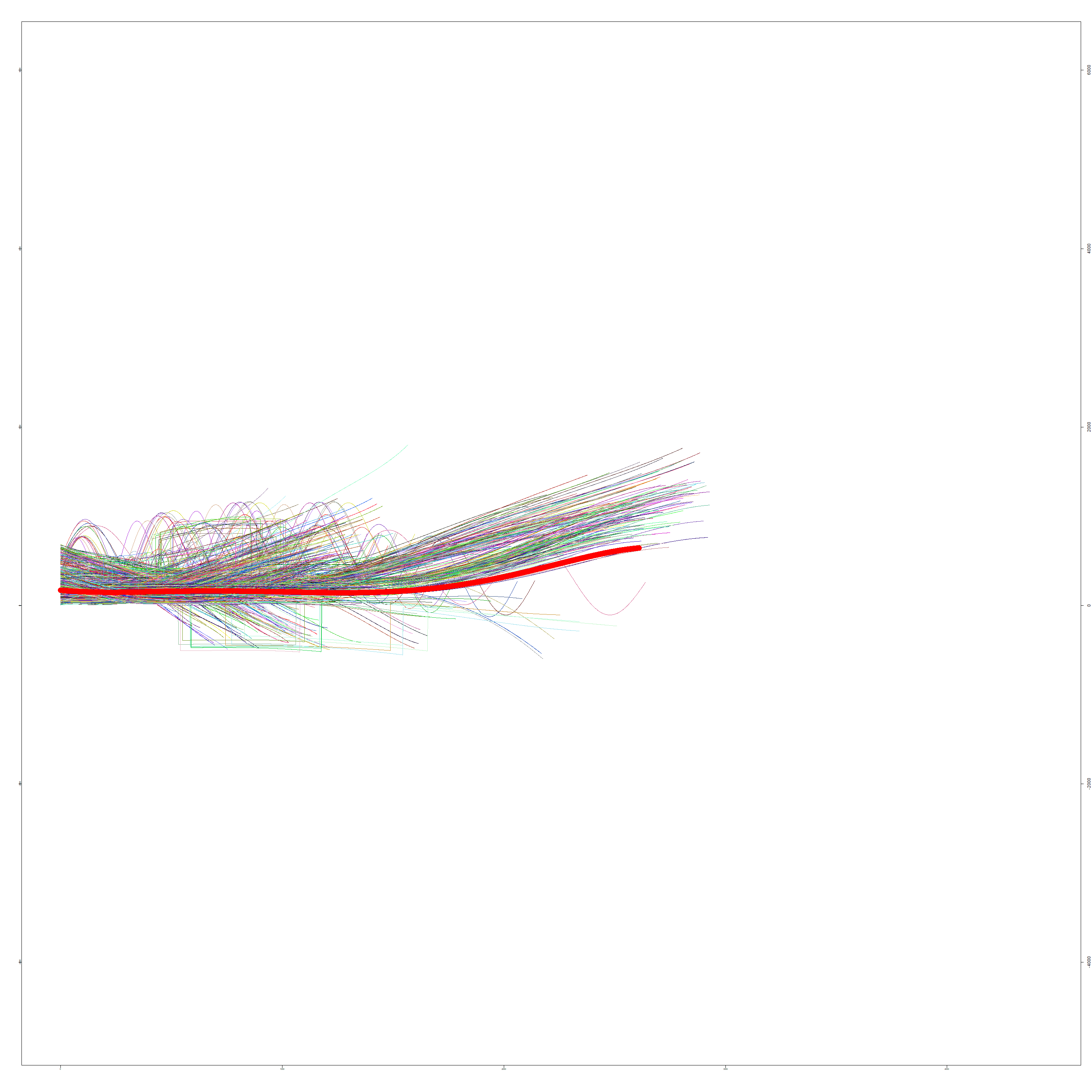} \\
$C_5$
\end{minipage} 
\\

\begin{minipage}[b]{0.19\textwidth}
\centering
\includegraphics[width=0.95\linewidth]{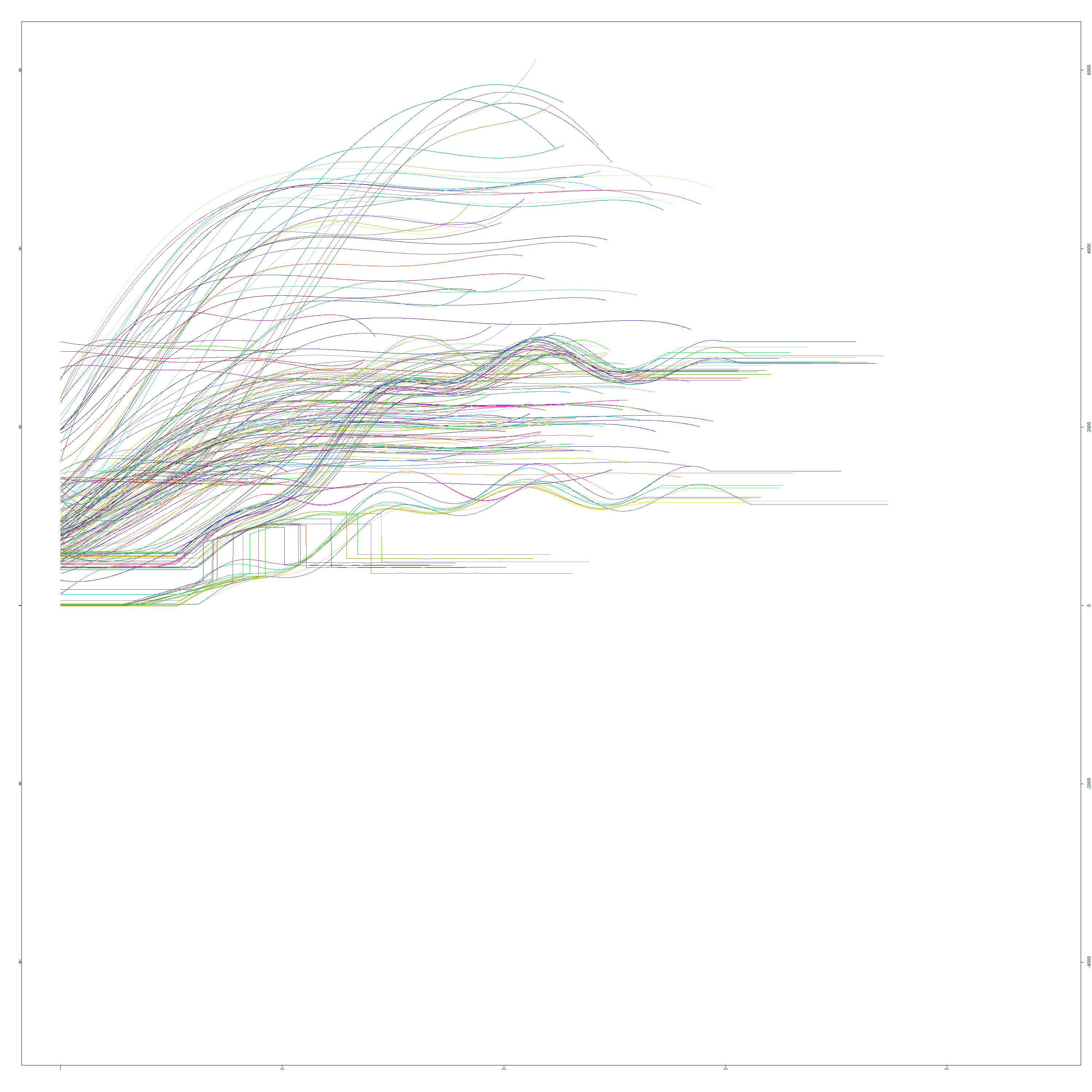} \\
$D_1$
\end{minipage} 
\begin{minipage}[b]{0.19\textwidth}
\centering
\includegraphics[width=0.95\linewidth]{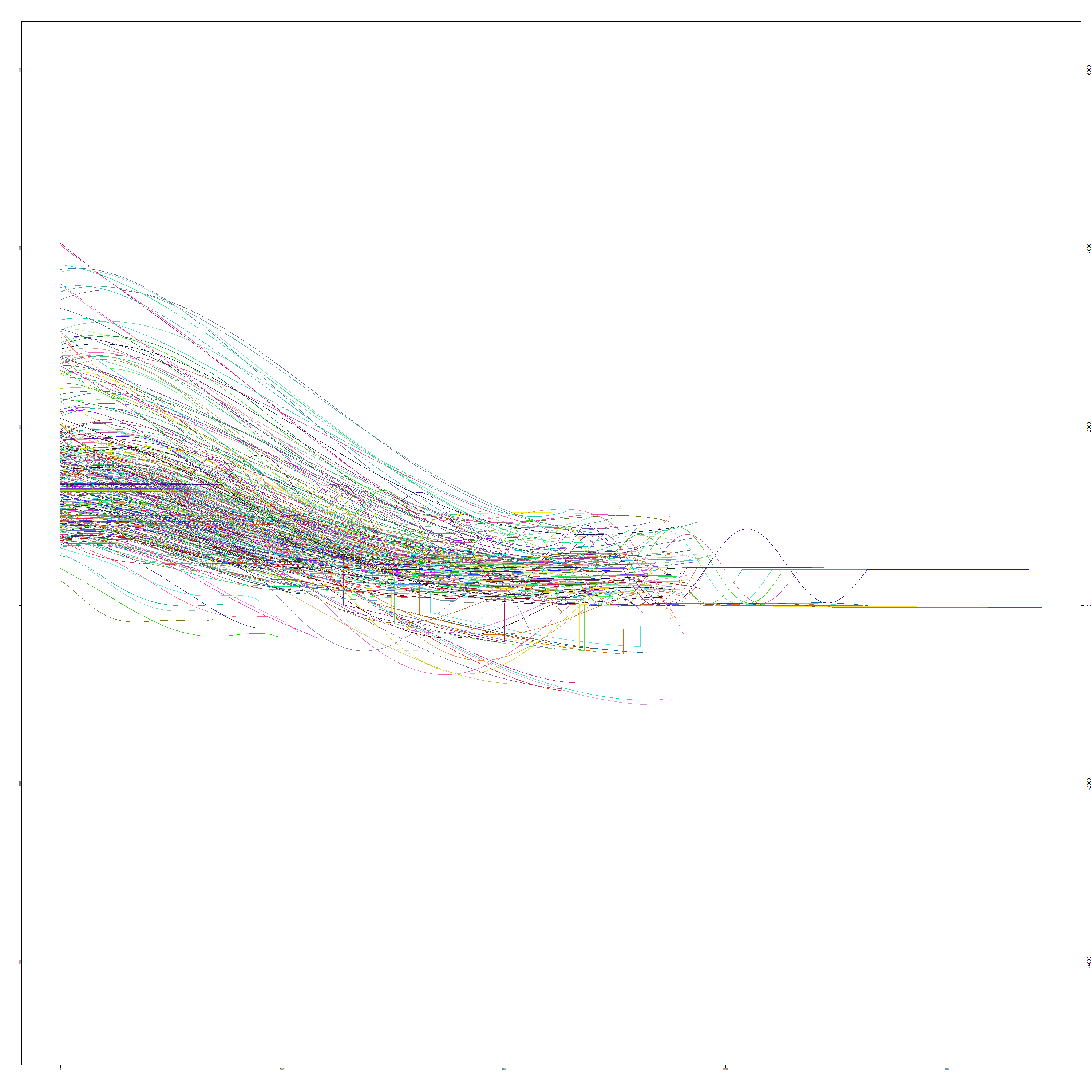} \\
$D_2$
\end{minipage}
\begin{minipage}[b]{0.19\textwidth}
\centering
\includegraphics[width=0.95\linewidth]{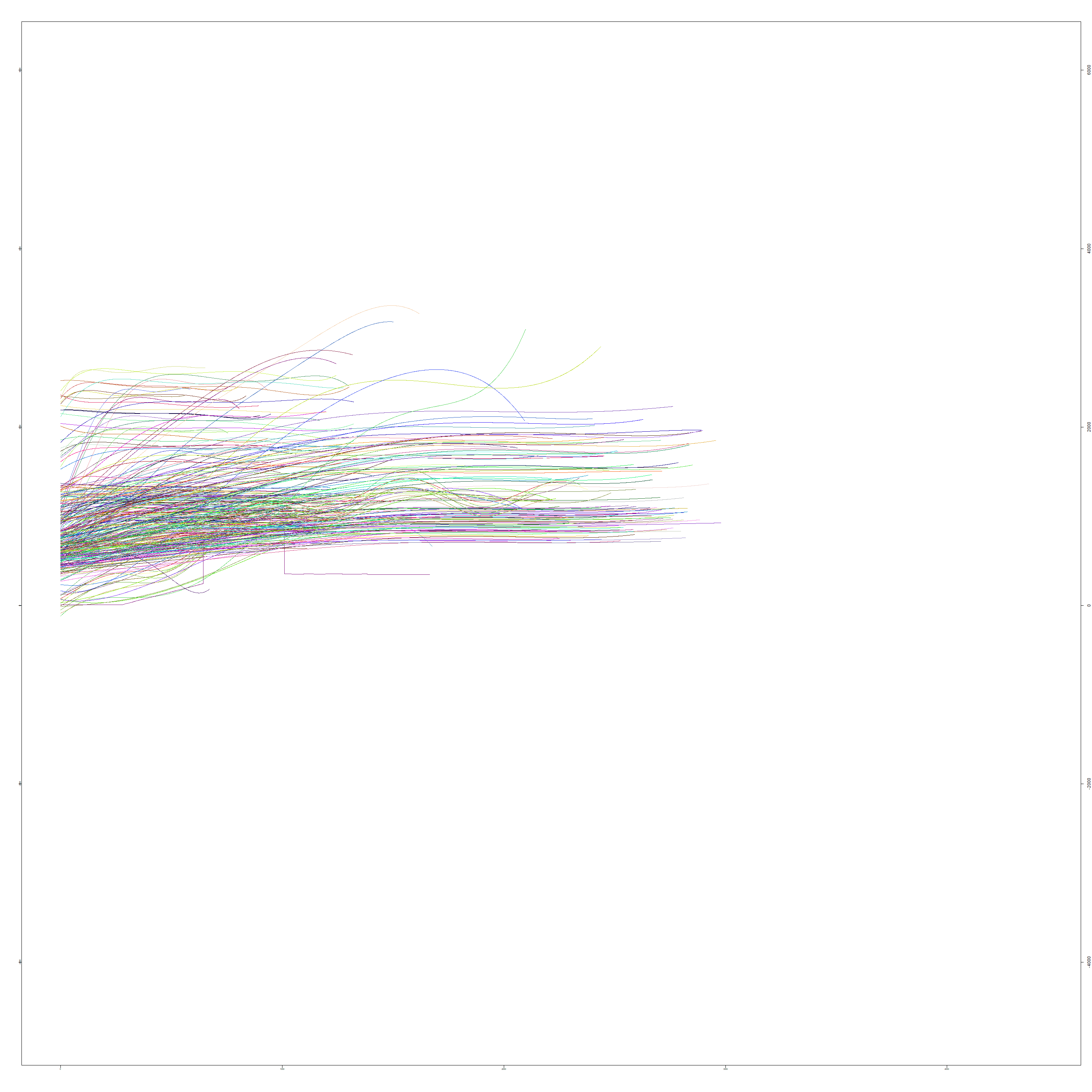} \\
$D_3$
\end{minipage} 
\begin{minipage}[b]{0.19\textwidth}
\centering
\includegraphics[width=0.95\linewidth]{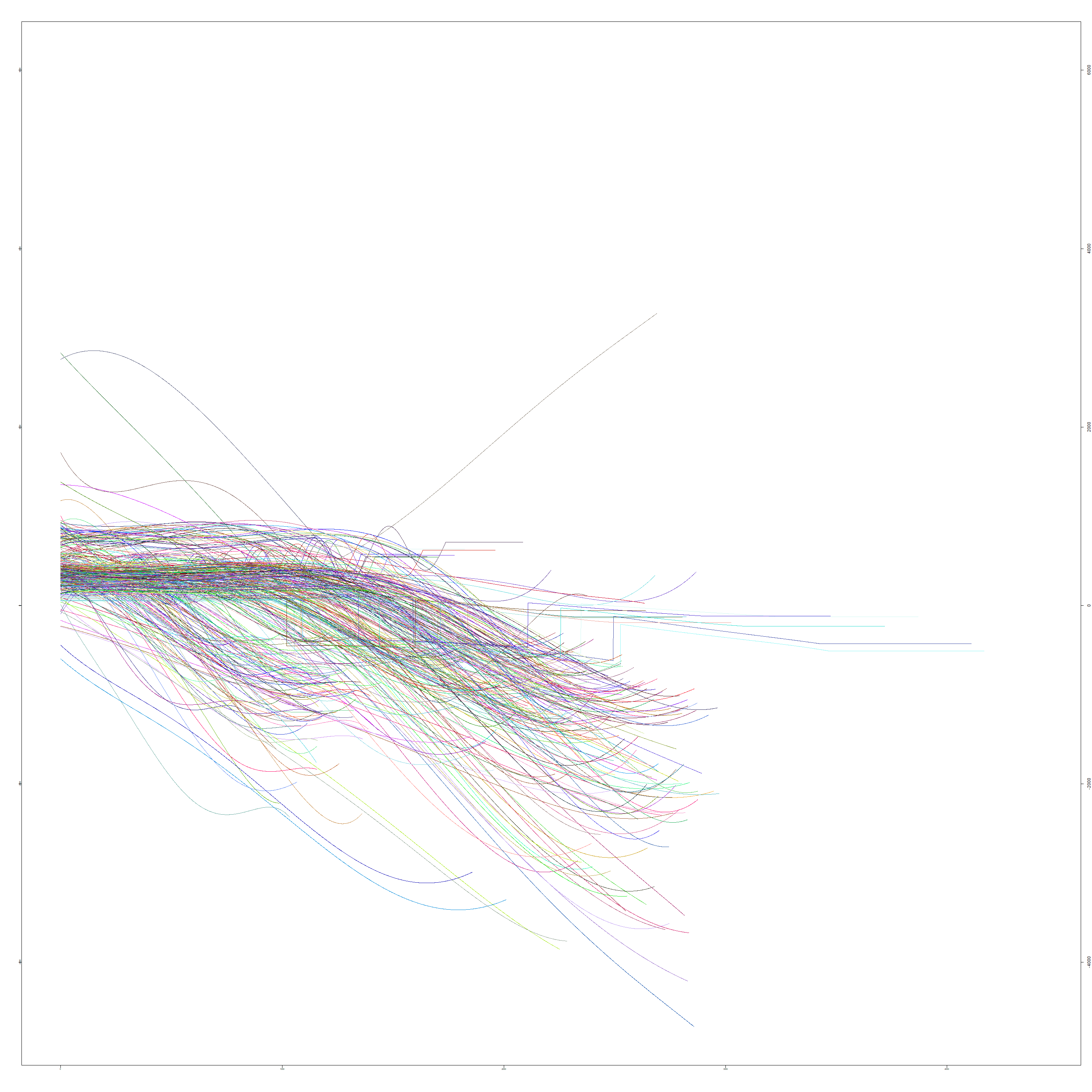} \\
$D_4$
\end{minipage} 
\begin{minipage}[b]{0.19\textwidth}
\centering
\includegraphics[width=0.95\linewidth]{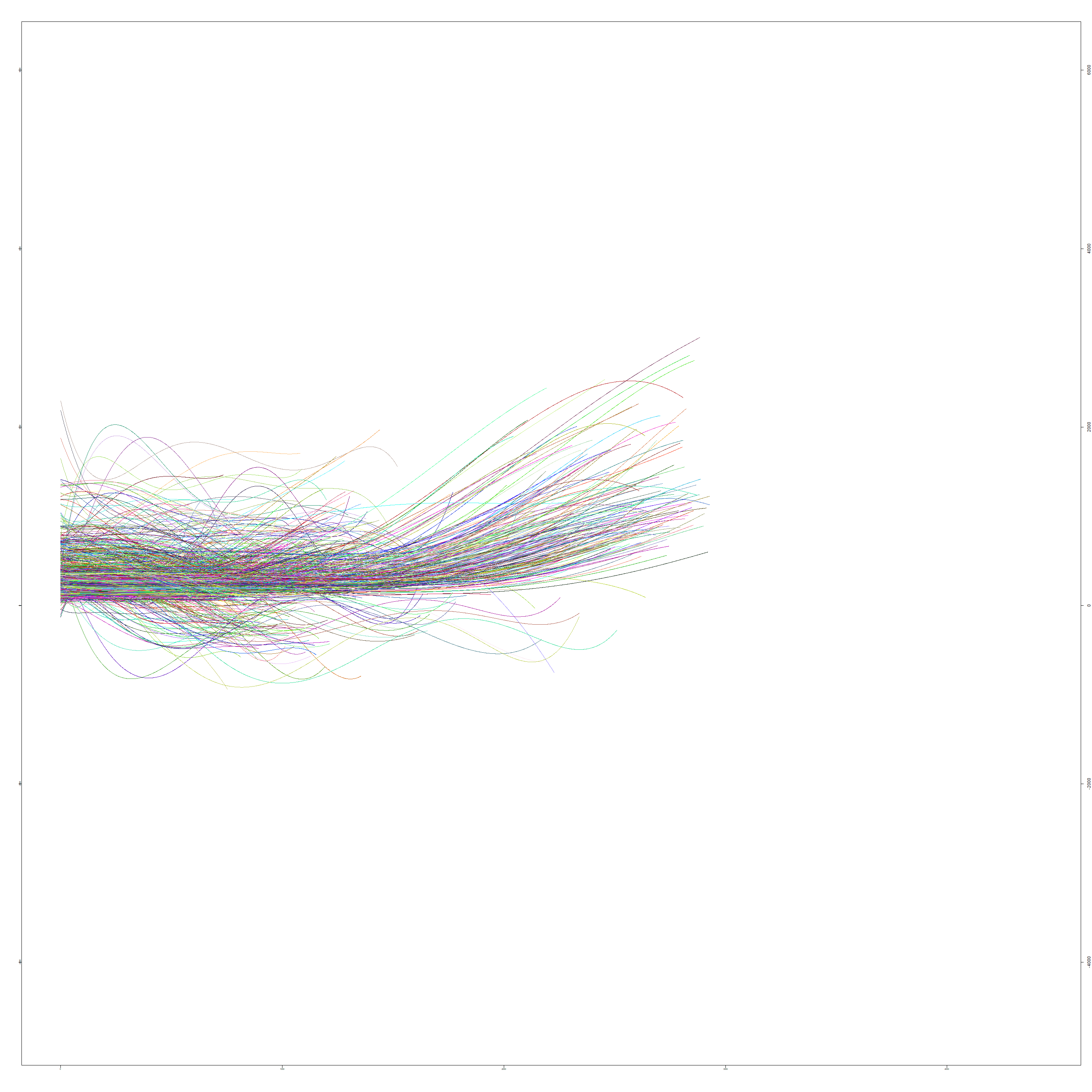} \\
$D_5$
\end{minipage} 
\caption{Clustering ($X$-curves on top and $Y$-curves below)} 
\label{Fig::all_clust}
\end{figure}

As explained in Section \ref{clustering}, we group the $X$-curves into five clusters denoted by $C_1, C_2, \ldots, C_5$ and displayed in Figure \ref{Fig::all_clust}. We observe that these clusters gather  curves with similar shapes and are homogeneous.

As for the $Y$-curves, they are grouped into sets $D_1, D_2, \ldots, D_5$ defined by
		$$D_i = \{Y_a / X_a \in C_i\}$$
for each $i= 1, \ldots, 5$. 


The next step consists in computing the reference curves $RC_1, RC_2, \ldots, RC_5$. They are drawn in red color in  Figure \ref{Fig::all_clust}.

\begin{figure}[H] 
\centering
\begin{minipage}[b]{0.20\textwidth}
\centering
\includegraphics[width=0.95\linewidth]{all_superclust_CR_1.png} \\
(a)
\end{minipage} 
\begin{minipage}[b]{0.20\textwidth}
\centering
\includegraphics[width=0.95\linewidth]{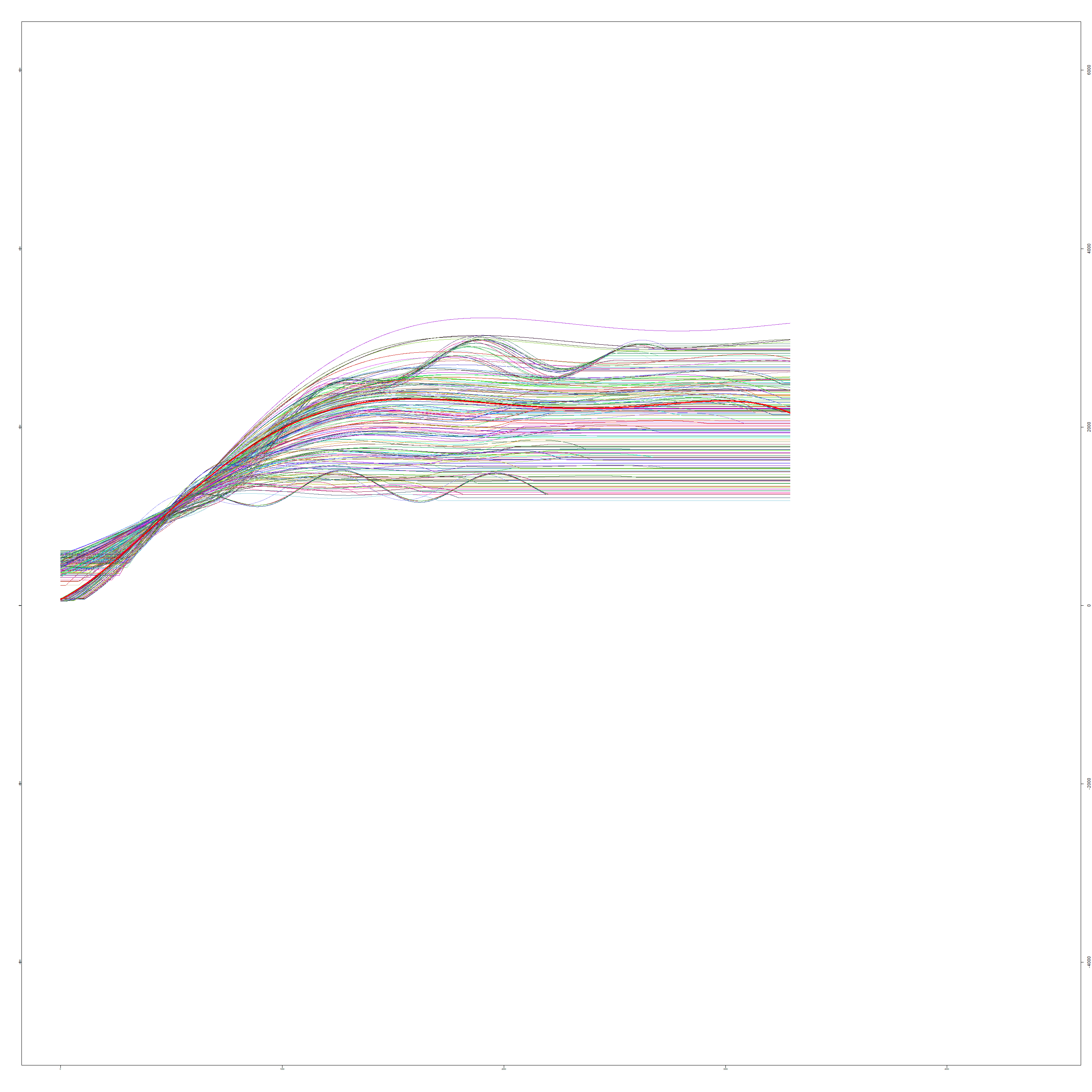} \\
(b)
\end{minipage} 
\begin{minipage}[b]{0.20\textwidth}
\centering
\includegraphics[width=0.95\linewidth]{bis_all_superclust_CR_1.png} \\
(c)
\end{minipage} 
\begin{minipage}[b]{0.20\textwidth}
\centering
\includegraphics[width=0.95\linewidth]{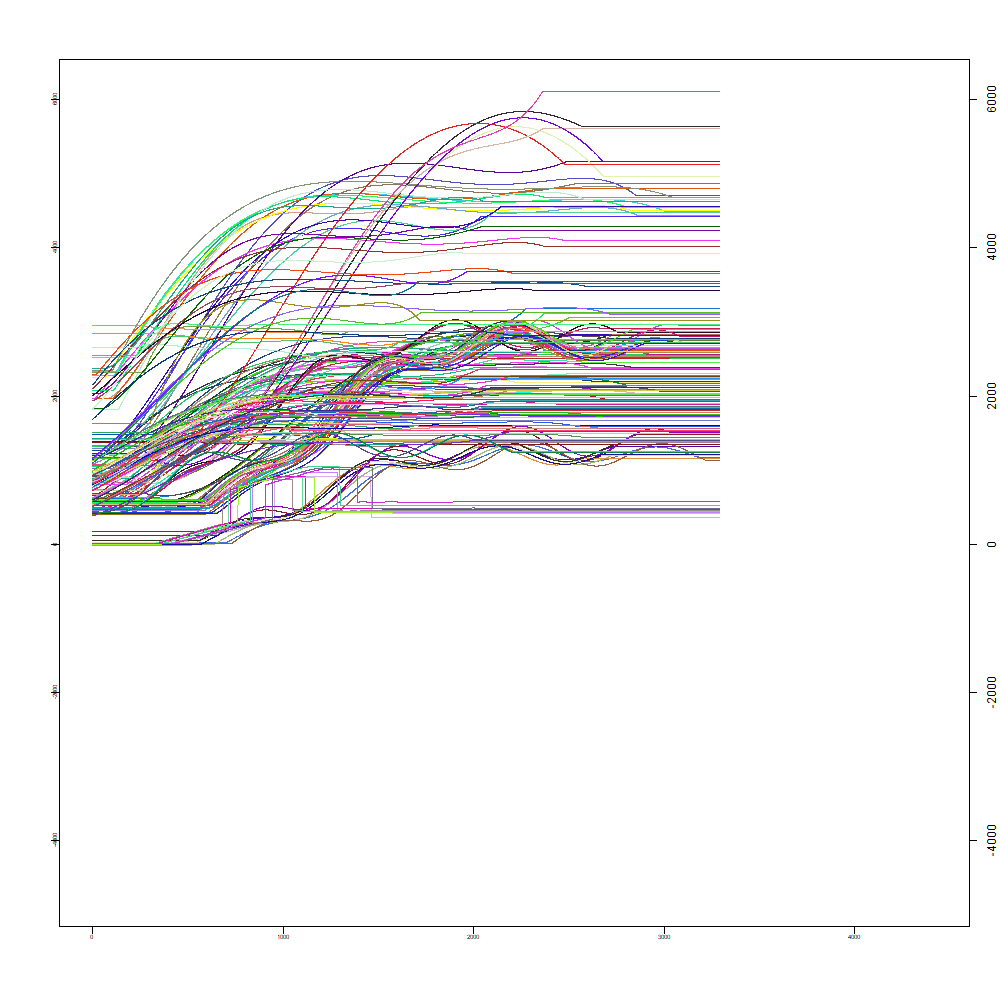} \\
(d)
\end{minipage}
\caption{Realignment and transformation of the curves in $C_1$ (a,b) and $D_1$ (c,d), initial curves in (a) and (c), transformed ones in (b) and (d)} 
\label{Fig::trans_D1}
\end{figure}

Then the realignment and the synchronization-transformation are  applied to the $X$- and $Y$- curves. To not be too long, the results of these transformations are displayed only for cluster $C_1$. Figures \ref{Fig::trans_D1}-a and \ref{Fig::trans_D1}-b contains the initial curves $X^i_a$ of cluster $C_1$ and their transformed $\breve{X}^1_a$.

In the same way, Figure \ref{Fig::trans_D1}-c and \ref{Fig::trans_D1}-d display the initial curves $Y^1_a$ of $D_1$ and their transformed $\breve{Y}^1_a$.
\vspace{-2em}
%

\vspace{-3em}

The last step consists in detecting the atypical curves in each cluster using the Confidence Tubes method (CT) and the Conditional Quantiles method (CQ) (see Section \ref{sec::CQ}). 

\begin{figure}[H] 
\centering
\begin{minipage}[b]{0.20\textwidth}
\centering
\includegraphics[width=0.95\linewidth]{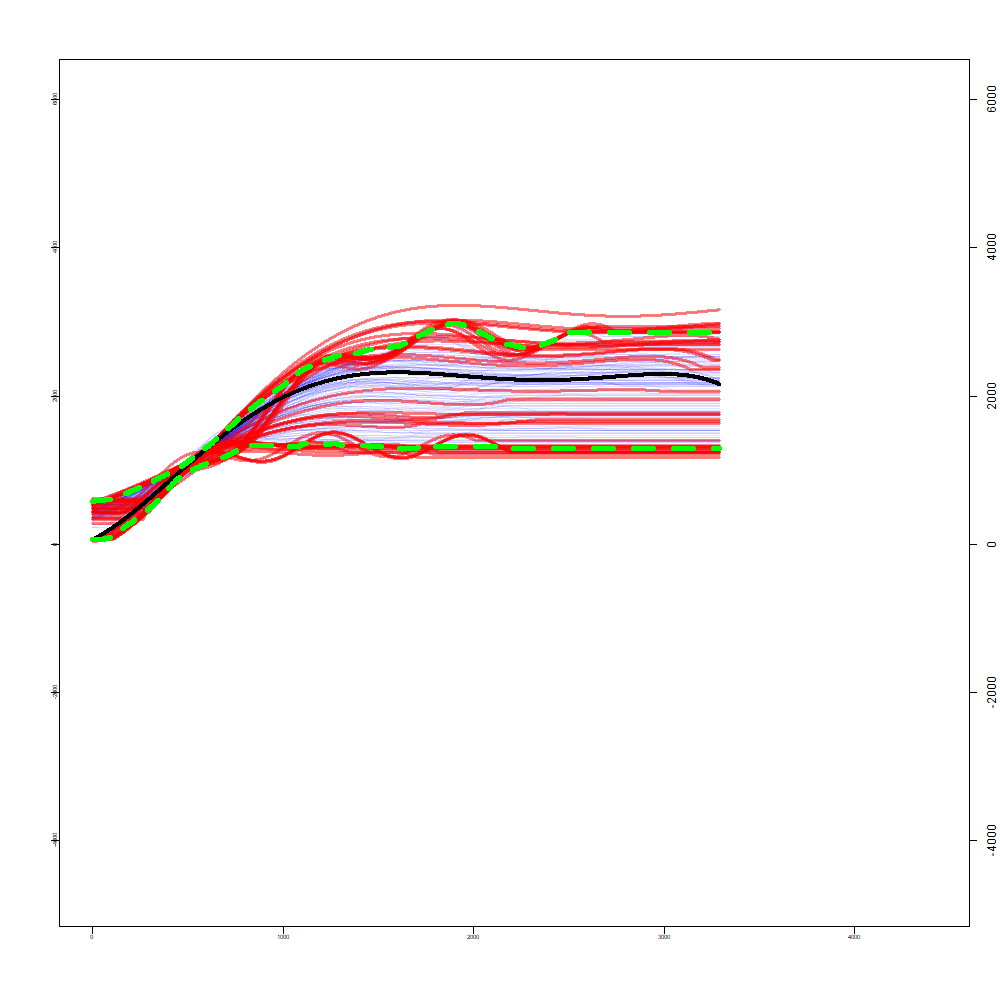} \\
(a)
\end{minipage} 
\begin{minipage}[b]{0.20\textwidth}
\centering
\includegraphics[width=0.95\linewidth]{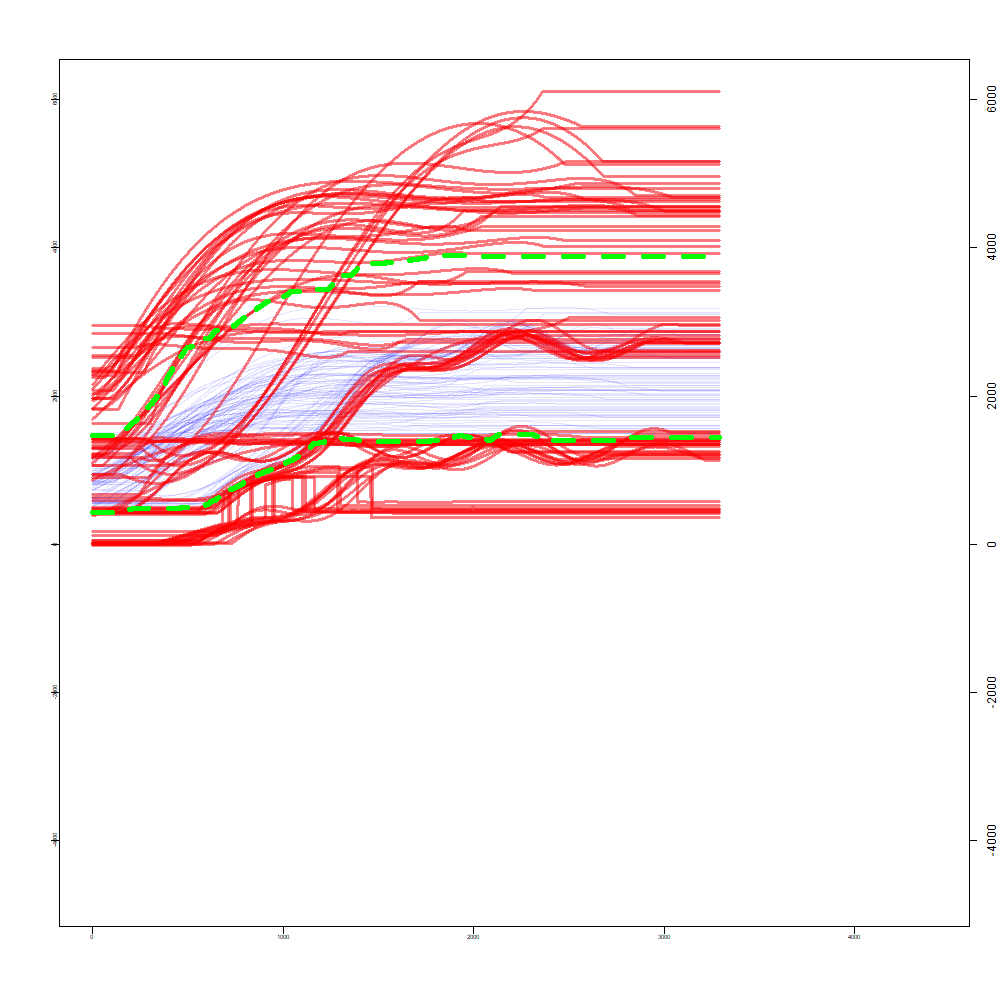} \\
(b)
\end{minipage} 
\begin{minipage}[b]{0.20\textwidth}
\centering
\includegraphics[width=0.95\linewidth]{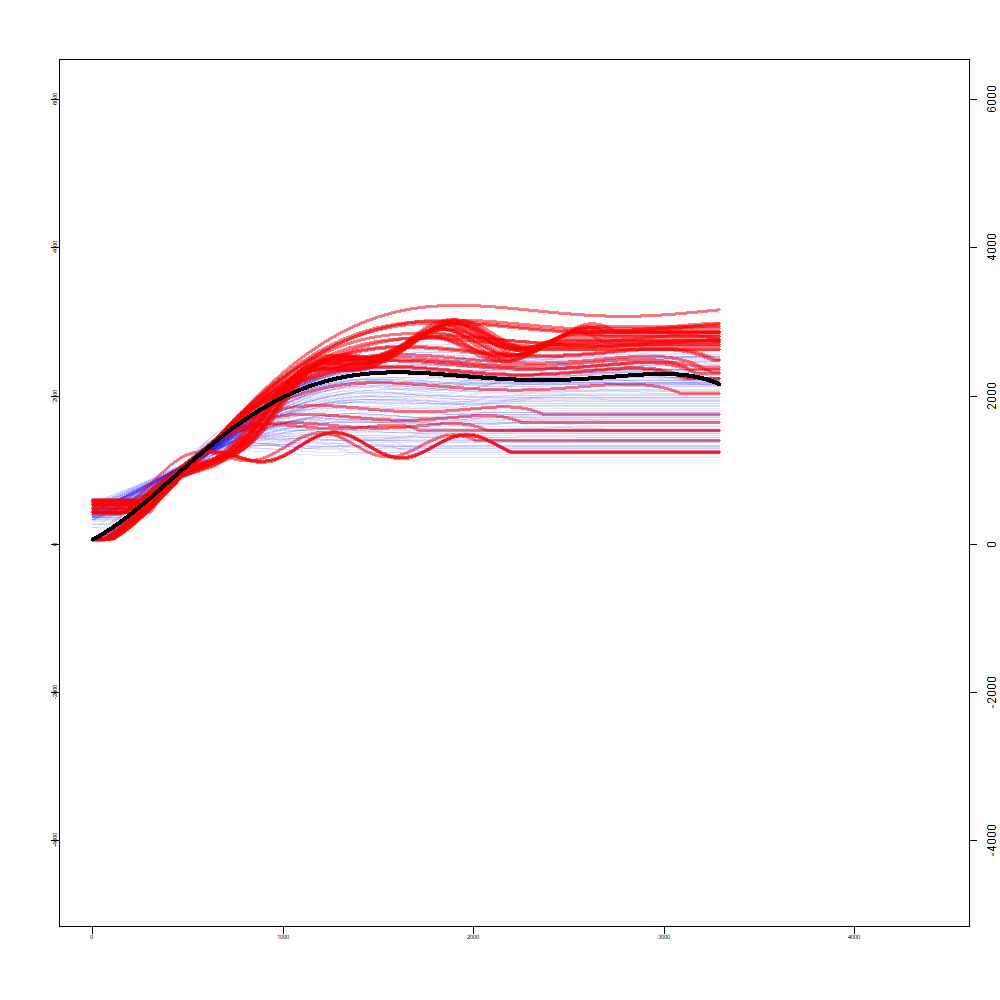} \\
(c)
\end{minipage} 
\begin{minipage}[b]{0.20\textwidth}
\centering
\includegraphics[width=0.95\linewidth]{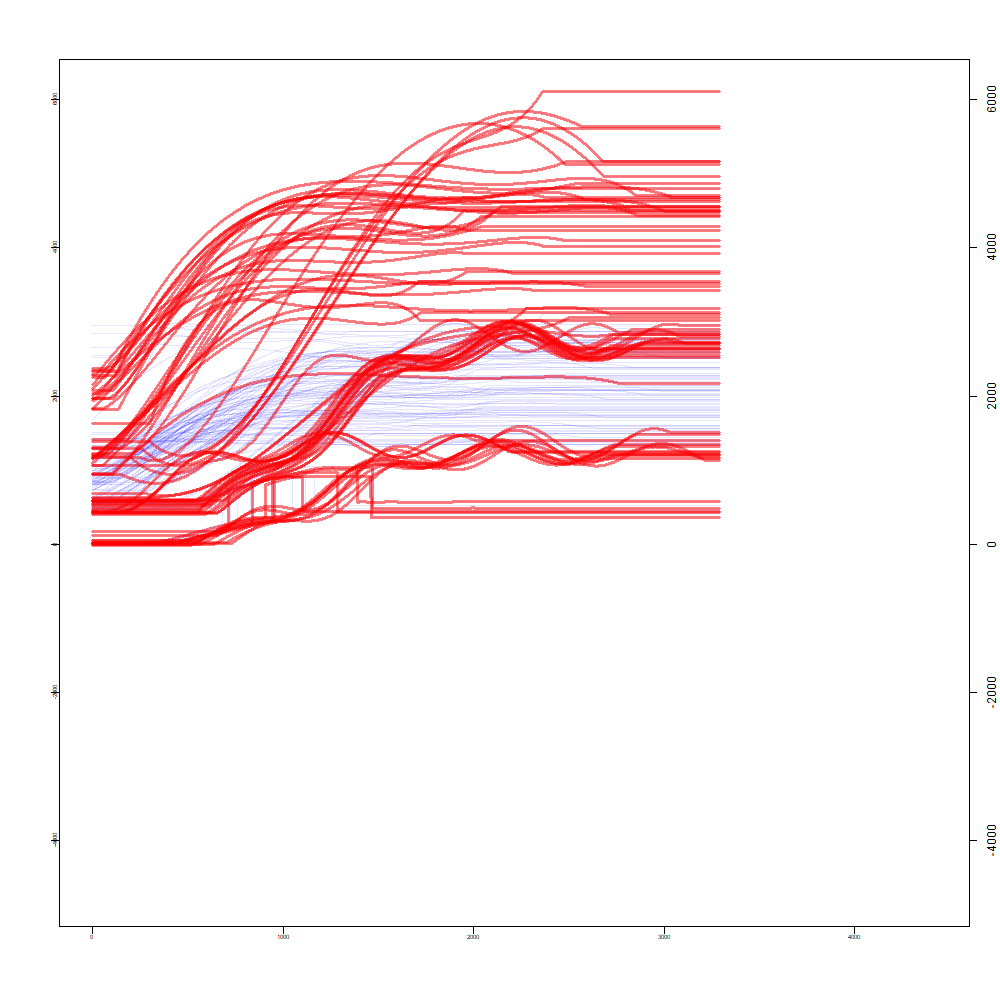} \\
(d)
\end{minipage}
\caption{All curves of $C_1$ and $D_1$ ($X$ on the left and $Y$ on the right) in blue and anomalies detected by CQ method in red} 
\label{Fig::anomali_C1_CQ}
\end{figure}

Figure \ref{Fig::anomali_C1_CQ}-a presents the results of the CT detection method for cluster $C_1$: all the $X$- curves in $C_1$ are drawn and the detected atypical ones are highlighted in red. The $Y$- curves and the detected atypical ones are displayed in Figure \ref{Fig::anomali_C1_CQ}-b. In the same way, Figure \ref{Fig::anomali_C1_CQ}-c and -d presents the result of the CQ detection method. At first glance, both methods seem to detect the same atypical curves !
\vspace{-2em}
\vspace{-2em}

\begin{table}[H]
\normalsize
\centering
\begin{tabular}{|c|c||c|c|c||c|c|c||c|c|c||c|c|c||c|c|c|}  
\hline
CT & $i$ & \multicolumn{3}{|c}{$1$} & \multicolumn{3}{|c}{$2$} & \multicolumn{3}{|c}{$3$} & \multicolumn{3}{|c}{$4$} & \multicolumn{3}{|c|}{$5$} \\ \hline
& & D & ND & T & D & ND & T & D & ND & T & D & ND & T & D & ND & T\\
\hline 
\hline 
$C_i$ & A & 15 & 9 & 24 & 22 & 3 & 25 & 26 & 0 & 26 & 19 & 0 & 19 & 41 &  25& 66\\ \hline 
  & NA & \textbf{27} & 131 & 158 & \textbf{20} & 297 & 317 & \textbf{11} & 272 & 283 & \textbf{11} & 291 & 302 & \textbf{30} &590 &620 \\ \hline 
  & T & 26 & 156 & 182 & 42 & 300 & 342  & 37 & 272 & 309 & 30 & 291 & 321 & 71 & 615 & 686\\ \hline 
  \hline 
  $D_i$ & A & 35 & 15 & 50 & 24 & 7 & 31 & 13 & 1 & 14 & 19 & 7 & 26 & 23 & 16 & 39 \\ \hline 
  & NA & \textbf{29} & 103 & 132 & \textbf{18} & 293 & 311  & \textbf{13} & 282 & 295 & \textbf{15} & 280 & 295 & \textbf{38} & 609 & 647\\ \hline 
  & T & 64 & 118 & 182 & 42 & 300 & 342  & 26 & 283 & 309 & 34 & 287 & 321 & 61 & 625 & 686\\ \hline 
  \hline 
    $C_i$ \& $D_i$ & A & 6 & 5 & 11 & 9 & 2 & 11 & 12 & 1 & 13 & 9 & 0 & 9 & 20 & 16 & 36\\ \hline 
  & NA & \textbf{11} & 160 & 171 & \textbf{10} & 321 & 331 & \textbf{10} & 286 & 296 & \textbf{10} & 302 & 312 & \textbf{16} & 634 & 650\\ \hline 
  & T & 17 & 165 & 182 & 19 & 323 & 342  & 22 & 287 & 309 & 19 & 302 & 321 & 36 & 650 & 686\\ \hline 
 
\end{tabular}
\caption{Confusion matrices for the CT method, in bold the number of false alarms} \label{tab:CT}
\end{table}

\begin{table}[H]
\normalsize
\centering
\begin{tabular}{|c|c||c|c|c||c|c|c||c|c|c||c|c|c||c|c|c|}  
\hline
CQ & $i$ & \multicolumn{3}{|c}{$1$} & \multicolumn{3}{|c}{$2$} & \multicolumn{3}{|c}{$3$} & \multicolumn{3}{|c}{$4$} & \multicolumn{3}{|c|}{$5$} \\ \hline
& & D & ND & T & D & ND & T & D & ND & T & D & ND & T & D & ND & T\\ \hline 
\hline 
$C_i$ & A  & 23 & 1 & 24 & 24 & 1 & 25 & 25 & 1 & 26 & 17 & 2 & 19 & 45 & 21 & 66\\ \hline 
      & NA & \textbf{0} & 158 & 158 & \textbf{1} & 316 & 317 & \textbf{5} & 278 & 283 & \textbf{5} & 297 & 302 & \textbf{25} & 595 & 620\\ \hline 
      & T  & 23 & 159 & 182 & 25 & 317 & 342 & 30 & 279 & 309 & 22 &  299& 321 & 70 & 616 & 686\\ \hline 
  \hline 
$D_i$ & A  & 37 & 13 & 50 & 26 & 5 & 31 & 14 & 0 & 14 & 23 & 3 & 26 & 26 & 13 & 39\\ \hline 
      & NA & \textbf{13} & 119 & 132 & \textbf{13} & 298 & 311 & \textbf{12} & 283 & 295 & \textbf{15} &  280& 295 & \textbf{35} & 612 & 647\\ \hline 
      & T  & 50 & 132 & 182 & 39 & 303 & 342 & 26 & 283 & 309 & 38 &  283 & 321 & 61 & 625 & 686\\ \hline 
  \hline 
$C_i$ \& $D_i$ & A  & 8 & 3 & 11 & 9 & 2 & 11 & 12 & 1 & 13 & 7 & 2 & 9 & 19 & 17 & 36\\ \hline 
      & NA & \textbf{7} & 164 & 171 & \textbf{8} & 323 & 331 & \textbf{8} & 288 & 296 & \textbf{10} & 302 & 312 & \textbf{19} & 631 & 650\\ \hline 
      & T  & 15 & 167 & 182 & 17 & 325 & 342 & 20 & 289 & 309 & 17 & 304 & 321 & 38 & 648 & 686\\ \hline

\end{tabular}
\caption{Confusion matrices for the CQ method, in bold the number of false alarms} \label{tab:CQ}
\end{table}

To display the confusion matrices for both methods, we use some abbreviate notations: A for atypical curves, NA for normal curves, D for  curves detected as atypical, ND for curves detected as normal. 

If we examine the first confusion matrices ( Table \ref{tab:CT}), which cross the true simulated anomalies and the detected ones by the CT  method, we see that there are a not negligible number of false alarms, for each cluster. They mainly correspond to normal curves having very high or very low values.

Table \ref{tab:CQ} shows the results that we get by using the CQ detection method. We observe that the results are much better, the number of false alarms decreases drastically. This fact can be explained because the CQ detection method takes into account the nature of the data which are time series, while the CT method considers the values as a simple cloud of points.



\section{Aircraft engine real data}

In this section, we deal with real bivariate data recorded  during the flights.  What follows is a part of the Cynthia Faure's PHD work, \cite{Faure2018}, made in collaboration with the Health Monitoring  Department of the Safran Aircraft Company of R\'eau, France.

In fact the sensors on the engine record more than 50 variables, but we took as an example the fan speed as \textit{key variable} $X$ and the temperature inside the motor as $Y$. The database contains the data relative to about 549 flights and 8 different engines. Each flight has a duration of around 2.8 hours.

4500 transient ascending phases were extracted from the 549 records. Their lengths are comprised between 200 and 10000 time units (8 Hz). We apply the methodology defined in Section 2, and seek to detect the anomalies as well for variable $X$ as for variable $Y$ or for both.

%
%
%
%
%
%
%

\begin{figure}[h!] 
\begin{minipage}[b]{0.24\linewidth}
\centering
\includegraphics[width=0.99\linewidth]{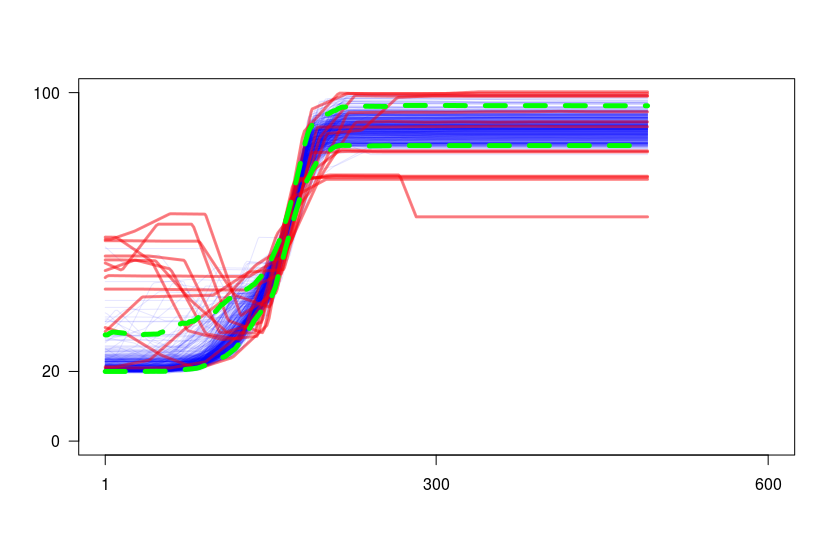}\\
(a) 
\end{minipage} 
\begin{minipage}[b]{0.24\linewidth}
\centering
\includegraphics[width=0.99\linewidth]{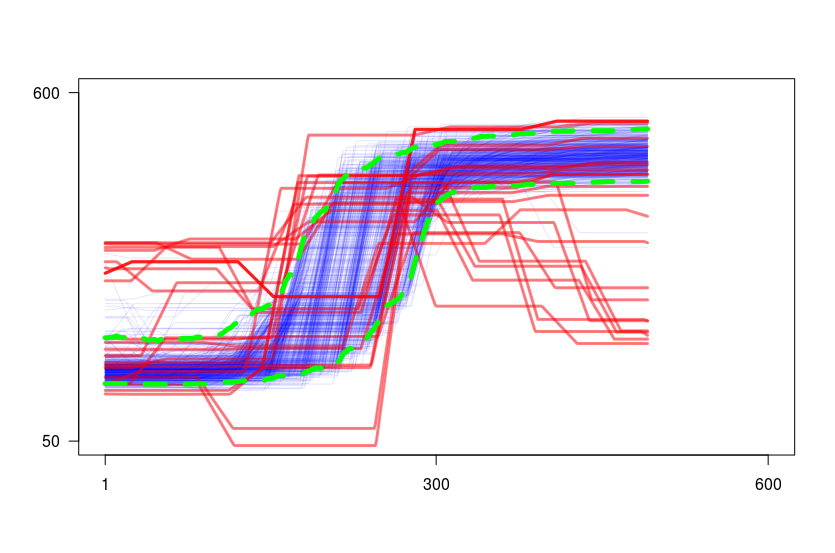} \\
(b)

\end{minipage} 
\begin{minipage}[b]{0.24\linewidth}
\centering
\includegraphics[width=0.95\linewidth]{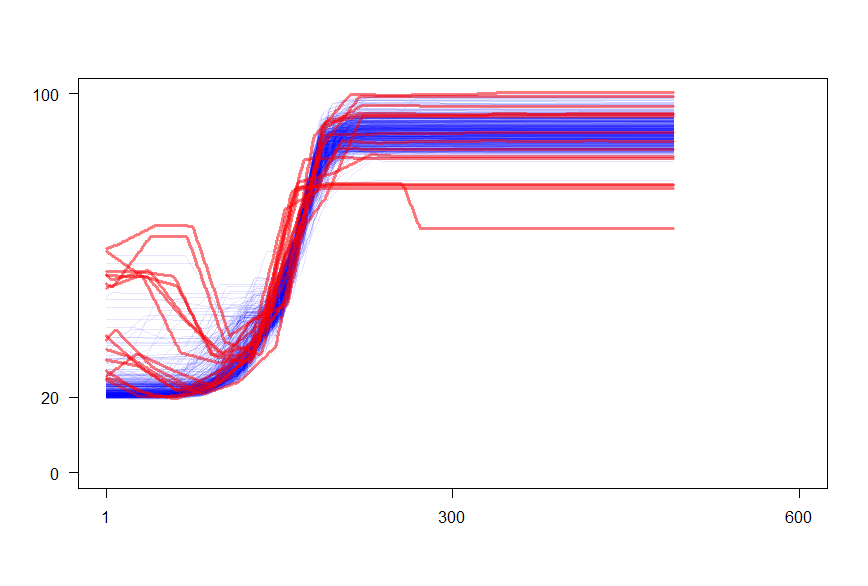}\\
(c)
 
\end{minipage}
\begin{minipage}[b]{0.24\linewidth}
\centering
\includegraphics[width=0.95\linewidth]{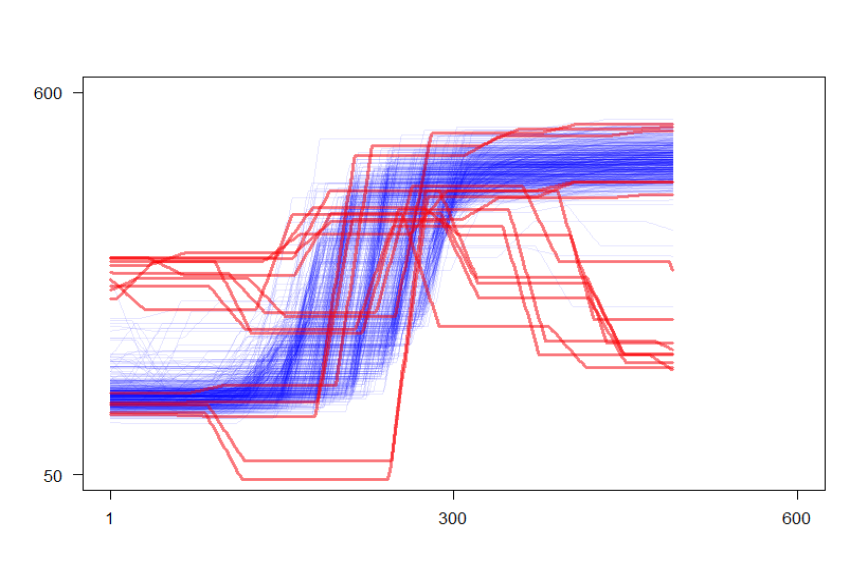}\\
(d)
\end{minipage} 

\caption{Results of the CT (a and b) and the CQ (c and d) methods: $X$-curves (a and c) and $Y$-curves  (b and d) are drawn,  the normal ones are in blue, the atypical ones in red} 
\label{Fig::N1_T3}
\end{figure}

Here for not be too long, we limit ourselves to one cluster, the one that mainly includes take-off phases. 
We then apply the CT and CQ anomaly detection methods.  In Figures \ref{Fig::N1_T3}), all the $X$- and $Y$-curves of this cluster are drawn in blue, except the detected atypical ones which are colored in red. At the top, the atypical curves are detected by the CT detection method, while at the bottom, they are detected by the CQ detection method.
 
%
%
%
%
%

As before, we see that both methods seem to detect nearly the same anomalies. In reality,  the CT method detects 12 atypical $X$-curves and CQ 14 ones (6 in common). As for the  $Y$-curves, CT detects 24 curves and CQ 18 curves of which 14 are common.

\begin{figure}[h!] 
\begin{minipage}[b]{0.24\linewidth}
\centering
\includegraphics[width=0.85\linewidth]{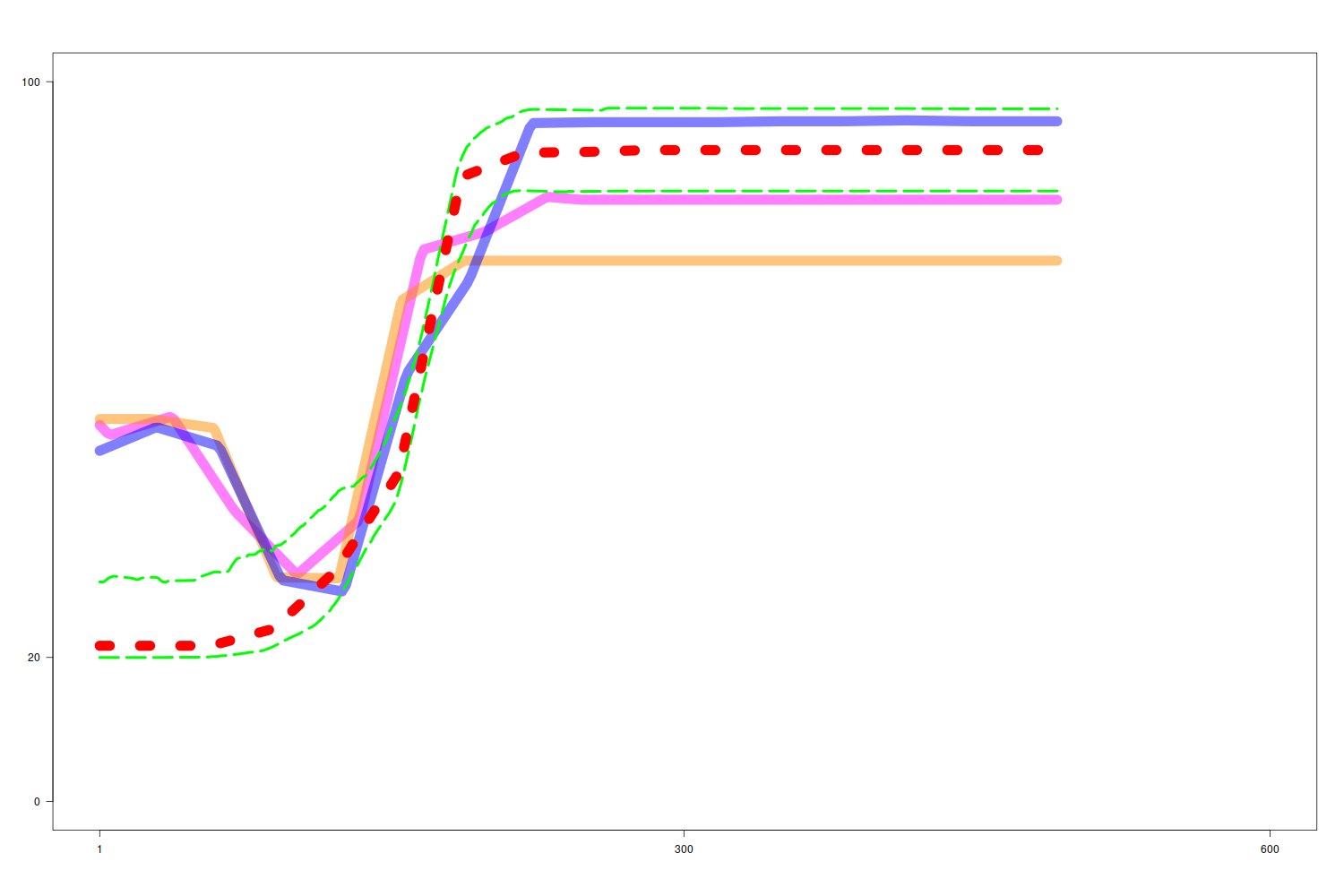} \\
(a) 
\end{minipage} 
\begin{minipage}[b]{0.24\linewidth}
\centering
\includegraphics[width=0.85\linewidth]{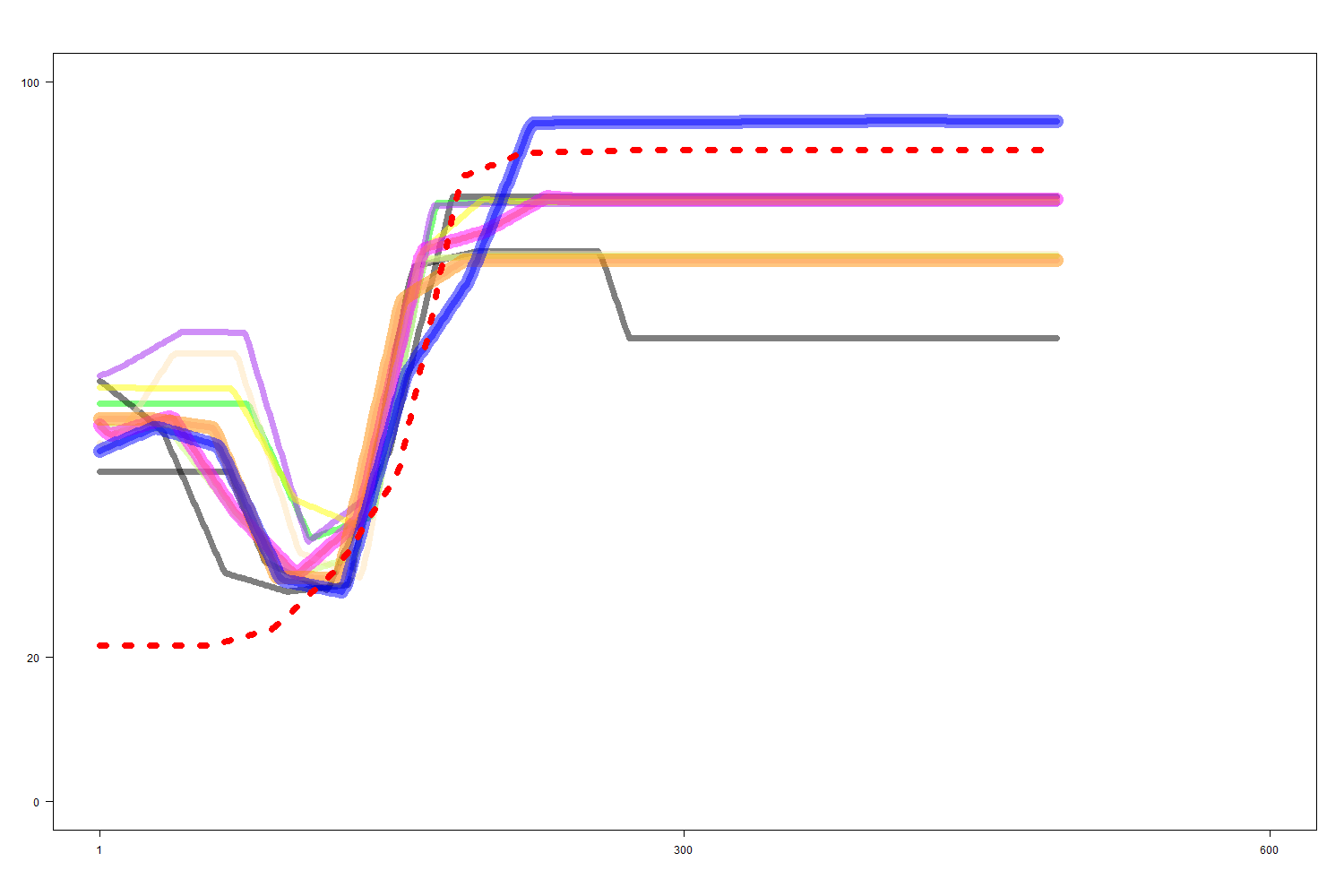} \\
(b) 
\end{minipage} 
\begin{minipage}[b]{0.24\linewidth}
\centering
\includegraphics[width=0.85\linewidth]{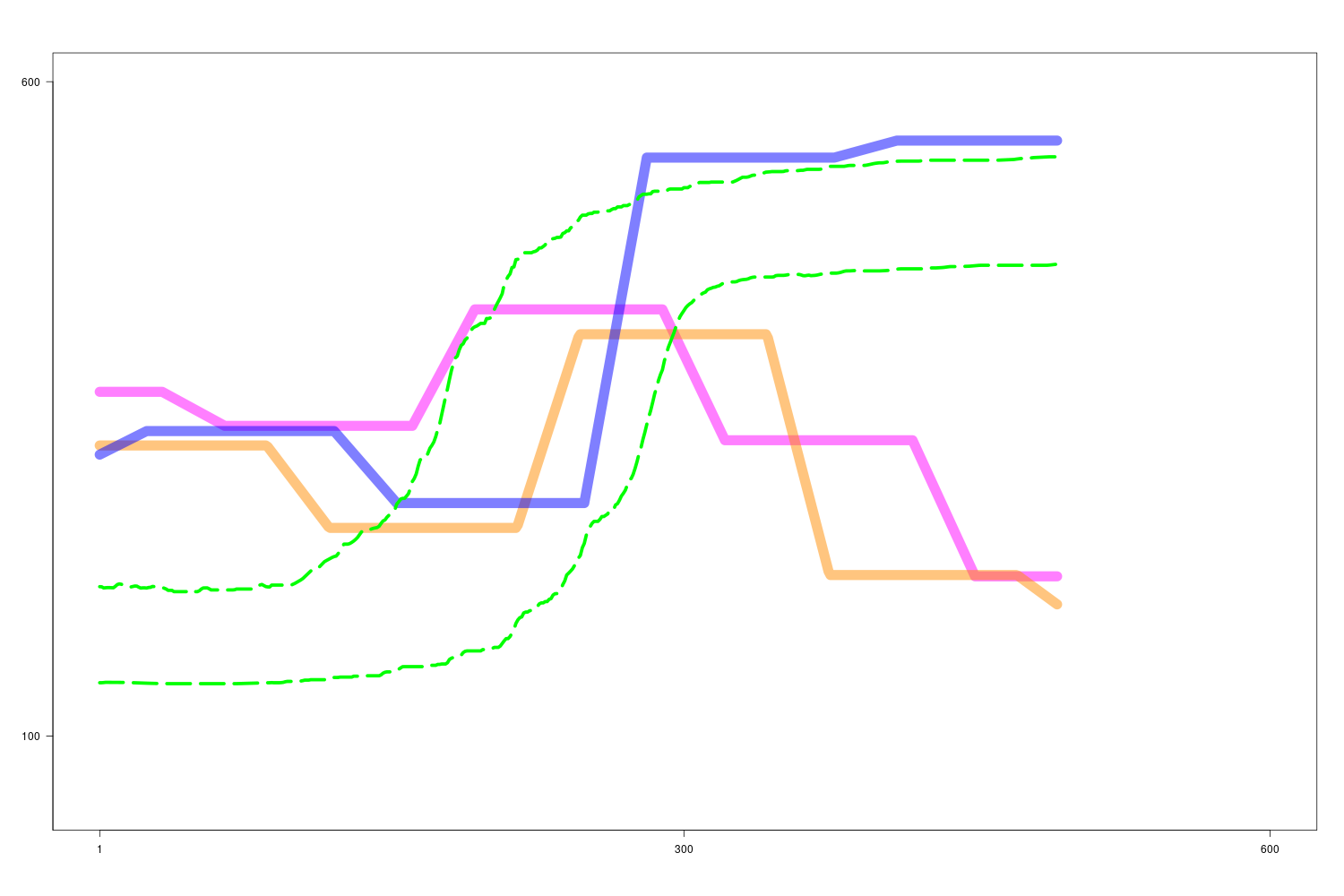}\\
(c) 
\end{minipage}
\begin{minipage}[b]{0.24\linewidth}
\centering
\includegraphics[width=0.85\linewidth]{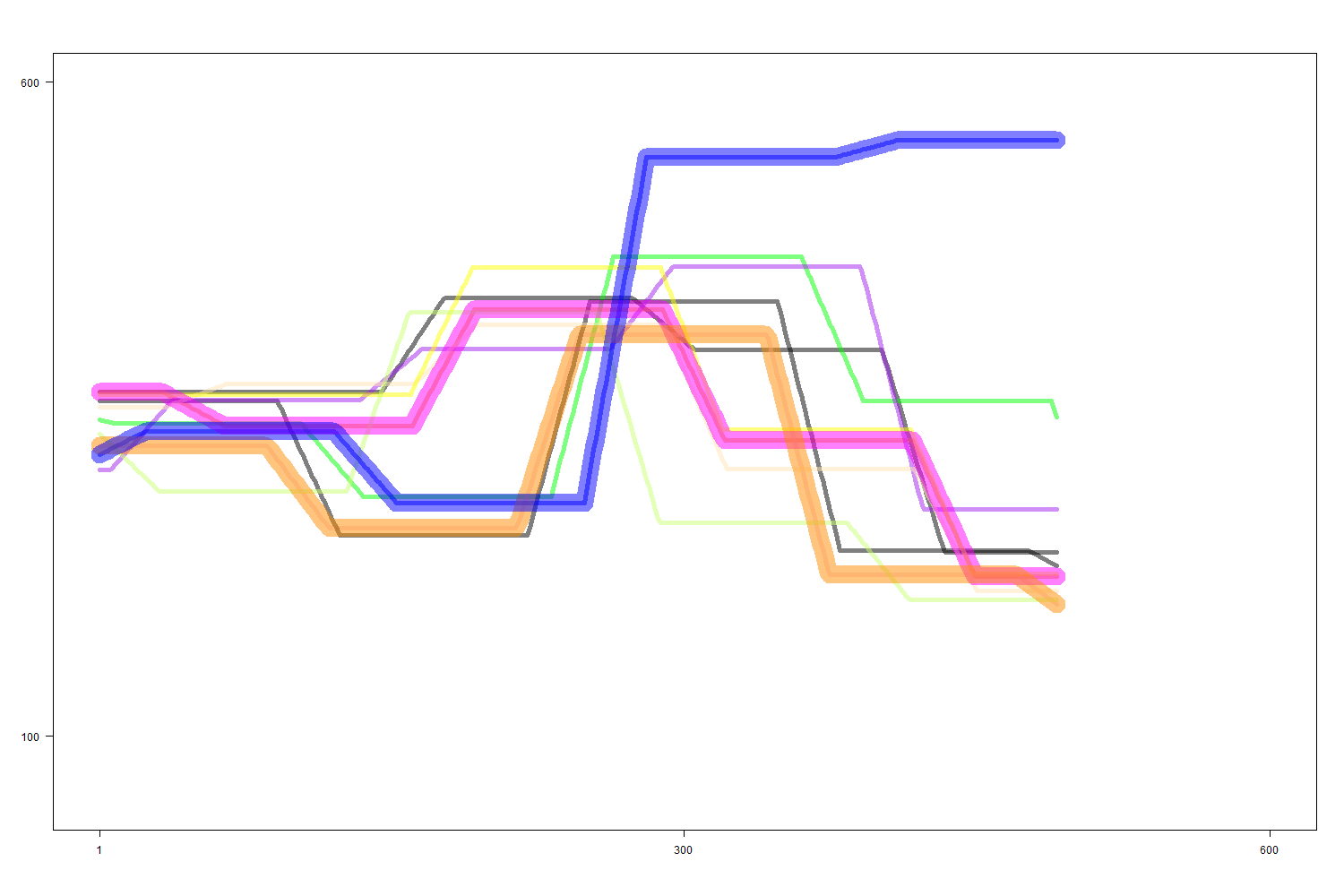}\\
(d) 
\end{minipage} 
\caption{Common atypical curves between X (a and b) and Y (c and d) with CT method (a and c) and CQ method (b and d). Reference curve in red and confidence tube in green for CT method.}
\label{Fig::N1_T3_common}
\end{figure}

Let us highlight the case of bidimensional curves which are detected as atypical for X \textbf{and}  for Y by both methods (see Figure \ref{Fig::N1_T3_common}).

The CT method detects three couples which are atypical for $X$ and for $Y$, while  the CQ method finds ten couples, that include the first three.

One more time, we can observe that the CQ anomaly detection has better performance than the CT method, even though it is not possible to compute the confusion matrices since we have no a priori knowledge about the existence of the anomalies. In this real-world study, the experts have been able to bring a validation to our findings. The detected cases corresponded to some events that they could identify.

\section{Conclusion}

This methodology provides very good results which help the experts to report abnormal behavior of important variables related to the engine operation. In the case of these aircraft engine real data, the companies can use these tools to increase  the probability of detecting any kind of atypical, abnormal behavior of some recorded variable, in order to prevent any incident and to plan the maintenance events.

\bibliography{biblio_iwann} 

\begin{thebibliography}{10}
\providecommand{\url}[1]{\texttt{#1}}
\providecommand{\urlprefix}{URL }

\bibitem{Baragona}
Baragona, R., Battaglia, F.: Outliers detection in multivariate time series by
  independent component analysis. Neural Computation  19(7),  1962--1984 (2007)

\bibitem{Bellas2014}
Bellas, A., Bouveyron, C., Cottrell, M., Lacaille, J.: Anomaly detection based
  on confidence intervals using som with an application to health monitoring.
  In: Villmann, T., Schleif, F.M., Kaden, M., Lange, M. (eds.) Advances in
  Self-Organizing Maps and Learning Vector Quantization, Proceedings of the
  10th International Workshop, (WSOM 2014). pp. 145--155. AISC,
  Springer-Verlag, Mittweida, Germany (July 2014)

\bibitem{Chandola2009}
Chandola, V., Banerjee, A., Kumar, V.: Anomaly detection: A survey. ACM Comput.
  Surv.  41(3),  15:1--15:58 (Jul 2009)

\bibitem{Charlier2015}
Charlier, I., Paindaveine, D., Saracco, J.: Conditional quantile estimation
  based on optimal quantization: from theory to practice. Computational
  Statistics and Data Analysis  91,  20--39 (2015),
  \url{https://hal.inria.fr/hal-01108504}

\bibitem{chen2008multi}
Chen, X.y., Zhan, Y.y.: Multi-scale anomaly detection algorithm based on
  infrequent pattern of time series. Journal of Computational and Applied
  Mathematics  214(1),  227--237 (2008)

\bibitem{Cheng}
{Cheng}, H., {Tan}, P., {Potter}, C., {Klooster}, S.: A robust graph-based
  algorithm for detection and characterization of anomalies in noisy
  multivariate time series. In: 2008 IEEE International Conference on Data
  Mining Workshops. pp. 349--358 (2008)

\bibitem{Faure2018}
Faure, C.: D\'etection de ruptures et identification des causes ou des
  sympt\^omes dans le fonctionnement des turbor\'eacteurs durant les vols et
  les essais. Ph.D. thesis, Universit\'e Paris 1 Panth\'eon-Sorbonne (2018)

\bibitem{Faure2017_b}
Faure, C., Bardet, J.M., Olteanu, M., Lacaille, J.: Design aircraft engine
  bivariate data phases using change-point detection method and self-organizing
  maps. In: Conference: ITISE - International work-conference on Time Series.
  University of Granada, Granada, Spain (September 2017)

\bibitem{Faure2017_a}
Faure, C., Bardet, J.M., Olteanu, M., Lacaille, J.: Using self-organizing maps
  for clustering anc labelling aircraft engine data phases. In: Ed., M.C. (ed.)
  12th International Workshop on Self-Organizing Maps and Learning Vector
  Quantization, Clustering and Data Visualization (WSOM+ 2017). pp. 1--8 (June
  2017)

\bibitem{Galeano}
Galeano, P., Peña, D., Tsay, R.S.: Outlier detection in multivariate time
  series by projection pursuit. Journal of the American Statistical Association
   101(474),  654--669 (2006)

\bibitem{Keogh2005}
{Keogh}, E., {Lin}, J., {Fu}, A.: Hot sax: efficiently finding the most unusual
  time series subsequence. In: Fifth IEEE International Conference on Data
  Mining (ICDM'05). pp. 226--233 (2005)

\bibitem{Leith}
{Knorn}, F., {Leith}, D.J.: Adaptive kalman filtering for anomaly detection in
  software appliances. In: IEEE INFOCOM Workshops 2008. pp. 1--6 (2008)

\bibitem{Lakhina}
Lakhina, A., Crovella, M., Diot, C.: Characterization of network-wide anomalies
  in traffic flows. In: Proceedings of the 4th ACM SIGCOMM Conference on
  Internet Measurement. pp. 201--206. IMC '04, ACM, New York, NY, USA (2004)

\bibitem{Lin2003}
Lin, J., Keogh, E., Lonardi, S., Chiu, B.: A symbolic representation of time
  series, with implications for streaming algorithms. In: Proceedings of the
  8th ACM SIGMOD Workshop on Research Issues in Data Mining and Knowledge
  Discovery. pp. 2--11. DMKD '03, ACM, New York, NY, USA (2003)

\bibitem{Malhotra2015}
Malhotra, P., Vig, L., Shroff, G., Agarwal, P.: Long short term memory networks
  for anomaly detection in time series. In: Verleysen, M. (ed.) European
  Symposium on Artificial Neural Networks, Computational Intelligence and
  Machine Learning (ESANN 2015). pp. 89--94. Bruges, Belgium (2015)

\bibitem{Michael00}
Michael, C.C., Ghosh, A.: Two state-based approaches to program-based anomaly
  detection. In: In Proceedings of the 16th Annual Computer Security
  Applications Conference. pp. 21--30. IEEE Computer Society (2000)

\bibitem{Olteanu2015}
Olteanu, M., Villa-Vialaneix, N.: On-line relational and multiple relational
  som. Neurocomputing  147(1),  15--30 (2015)

\bibitem{Rabenoro2014}
Rabenoro, T., Lacaille, J., Cottrell, M., Rossi, F.: Anomaly detection based on
  indicators aggregation. In: International Joint Conference on Neural Networks
  (IJCNN 2014). pp. 2548--2555. Beijing, China (July 2014)

\bibitem{Samanta1989}
Samanta, T.: Non-parametric estimation of conditional quantiles. Statistics and
  Probability Letters  7,  497--412 (1989)

\end{thebibliography}
\bibliographystyle{splncs04}
\end{document}